\documentclass[12pt]{article}
\usepackage{amsmath, amsfonts, amssymb,amsthm, amscd}

\usepackage{pictexwd,dcpic}
\usepackage{graphicx}
\usepackage{epsf}
\usepackage{epsfig}
\usepackage{rlepsf}
\usepackage{psfrag}
\usepackage{color}
\usepackage{hyperref}
\usepackage{makeidx}

\theoremstyle{definition}

\theoremstyle{remark}

\newtheoremstyle{thm}% name
  {12pt}%      Space above, empty = `usual value'
  {12pt}%      Space below
  {\itshape}% Body font
  {\parindent}%         Indent amount (empty = no indent, \parindent = para indent)
  {\scshape}% Thm head font
  {.}%        Punctuation after thm head
  {5pt}% Space after thm head: \newline = linebreak
  {}%         Thm head spec

\theoremstyle{thm}
\newtheorem*{T6.17}{Theorem \ref{gradientthm1}}
\newtheorem{thm}{Theorem}[section]

\newtheoremstyle{prop}% name
  {12pt}%      Space above, empty = `usual value'
  {12pt}%      Space below
  {\itshape}% Body font
  {\parindent}%         Indent amount (empty = no indent, \parindent = para indent)
  {\scshape}% Thm head font
  {.}%        Punctuation after thm head
  {5pt}% Space after thm head: \newline = linebreak
  {}%         Thm head spec

\theoremstyle{prop}
\newtheorem{prop}[thm]{Proposition}
\newtheoremstyle{lem}% name
  {12pt}%      Space above, empty = `usual value'
  {12pt}%      Space below
  {\itshape}% Body font
  {\parindent}%         Indent amount (empty = no indent, \parindent = para indent)
  {\scshape}% Thm head font
  {.}%        Punctuation after thm head
  {5pt}% Space after thm head: \newline = linebreak
  {}%         Thm head spec

\theoremstyle{lem}
\newtheorem*{L2.6}{Lemma \ref{lem2.4}}
\newtheorem{lem}[thm]{Lemma}

\newtheoremstyle{defn}% name
  {12pt}%      Space above, empty = `usual value'
  {12pt}%      Space below
  {\itshape}% Body font
  {\parindent}%         Indent amount (empty = no indent, \parindent = para indent)
  {\scshape}% Thm head font
  {.}%        Punctuation after thm head
  {5pt}% Space after thm head: \newline = linebreak
  {}%         Thm head spec

\theoremstyle{defn}
\newtheorem{defn}[thm]{Definition}

\newtheoremstyle{examp}% name
  {12pt}%      Space above, empty = `usual value'
  {12pt}%      Space below
  {}% Body font
   {\parindent}%         Indent amount (empty = no indent, \parindent = para indent)
  {\scshape}% Thm head font
  {.}%        Punctuation after thm head
  {5pt}% Space after thm head: \newline = linebreak
  {}%         Thm head spec

\theoremstyle{examp}
\newtheorem{examp}[thm]{Example}

\newtheoremstyle{cor}% name
  {12pt}%      Space above, empty = `usual value'
  {12pt}%      Space below
  {\itshape}% Body font
  {\parindent}%         Indent amount (empty = no indent, \parindent = para indent)
  {\scshape}% Thm head font
  {.}%        Punctuation after thm head
  {5pt}% Space after thm head: \newline = linebreak
  {}%         Thm head spec

\theoremstyle{cor}

\newtheoremstyle{recipe}% name
  {12pt}%      Space above, empty = `usual value'
  {12pt}%      Space below
  {\itshape}% Body font
   {\parindent}%         Indent amount (empty = no indent, \parindent = para indent)
  {\scshape}% Thm head font
  {.}%        Punctuation after thm head
  {5pt}% Space after thm head: \newline = linebreak
  {}%         Thm head spec

\theoremstyle{recipe}

\newtheoremstyle{rem}% name
  {12pt}%      Space above, empty = `usual value'
  {12pt}%      Space below
  {}% Body font
   {\parindent}%         Indent amount (empty = no indent, \parindent = para indent)
  {\scshape}% Thm head font
  {.}%        Punctuation after thm head
  {5pt}% Space after thm head: \newline = linebreak
  {}%         Thm head spec

\theoremstyle{rem}
\newtheorem{rem}[thm]{Remark}

\newcommand{\bp}{\begin{proof}}
\newcommand{\ep}{\end{proof}}

\newcommand{\norm}[1]{\left\Vert#1\right\Vert}
\newcommand{\abs}[1]{\left\vert#1\right\vert}

\newcommand{\ssc}{\text{sc}}

\renewcommand{\epsilon}{\varepsilon}

\newcommand{\id}{\operatorname{id}}
\newcommand{\ind}{\operatorname{Ind}}

\newcommand{\supp}{\operatorname{supp}}

\providecommand{\ker}[1]{$\text{ker}\ {#1}$}

\newcommand{\Q}{{\mathbb Q}}
\def\abs#1{\mathopen|#1\mathclose|}

\def\norm#1{\mathopen\|#1\mathclose\|}

{\catcode`"=12 \gdef\hex{"}}

\mathchardef\laplace=\hex0001

\mathchardef\nabla=\hex0272

\makeatletter

\def\@@dalembert#1#2{\setbox0\hbox{$#1\mathrm I$}

  \vrule height\ht0 depth\z@ width.04\ht0

  \rlap{\vrule height\ht0 depth-.96\ht0 width.8\ht0}

  \vrule height.1\ht0 depth\z@ width.8\ht0

  \vrule height\ht0 depth\z@ width.1\ht0 }

\def\dalembert{\mathbin{\mathpalette\@@dalembert{}}\,}

\makeatother

\begin{document}

%\frontmatter
\title{Polyfolds And A General Fredholm Theory}

\author{ H. Hofer\footnote{Research partially supported
by NSF grant  DMS-0603957.}
}

\maketitle
\tableofcontents

\section{Introduction}

In this paper we discuss the generalized Fredholm theory in polyfolds
as developed in \cite{HWZ1,HWZ2,HWZ3}\footnote{Here we shall not discuss the paper \cite{HWZ4}, which is the most relevant for symplectic field theory.}. Some of the results have been previously described in \cite{H1}. We  illustrate the concepts by the Gromov-Witten invariants and  refer the reader
to the upcoming papers \cite{HWZ6,HWZ7} and the lecture notes \cite{H2} and \cite{H3} for more details\footnote{We, i.e. HWZ, ultimately hope to finish our books
\cite{HWZ-polyfolds1} and \cite{HWZ-polyfolds2}, which will explain the theory
via many examples, though the key-material will already be contained
in the before-mentioned papers and lecture-notes.}.

The main point of this paper is to familiarize the reader with the language and the abstract results.
The theory has been applied to Gromov-Witten, Floer-theory, and symplectic field theory (SFT) and these applications are described
in the  above papers and lecture notes, but it is clear that the theory is applicable to other problems as well.
The theory has been devised to package moduli problems with analytical limiting behaviors, like bubbling-off, breaking of trajectories and stretching the neck, and to put them
into an abstract framework, which provides a Sard-Smale perturbation theory and deals with the transversality issues. An important point
is that the problems of interest very often exhibit compactness issues, but on the other hand allow very elaborate compactifications, which are the source for
interesting algebraic invariants, see for example \cite{G,BEHWZ}. In this paper we shall not discuss the Fredholm theory with operations which is needed for applications to Floer-theory and SFT. This theory will be described in the upcoming paper \cite{HWZ4} and the lecture notes \cite{H3}.

The aforementioned analytical limiting  phenomena, even assuming a sufficient amount of genericity, do not look like smooth phenomena if smoothness refers to the usual concept.
However, it turns out that the notion of smoothness can be relaxed, and a generalization of differential geometry and functional analysis can be developed, so that the limiting phenomena can be viewed as smooth phenomena, even if they are quite often obscured by transversality issues.
In this generalized context the classical nonlinear Fredholm theory can be extended to a much larger class of spaces and operators, which can deal with the aforementioned problems.

The starting point for our considerations is the  observation that the notion of differentiability in finite dimensions, which usually is generalized as Fr\'echet differentiability to infinite-dimensional Banach spaces,  can be generalized in a quite  different way if  Banach spaces are equipped with an additional piece of data.  This additional piece of structure, called an sc-structure,
occurs in interpolation theory, \cite{Tr}, albeit under a different name. In fact, we give a quite different interpretation of such a structure and make clear, that it can be viewed
as a generalization of a smooth structure on a Banach space. We call this generalization sc-smoothness, and the generalization of differentiability of a map we refer to as sc-differentiability.

The interesting thing is then the following fact. There are many sc-differentiable maps $r:U\rightarrow U$ satisfying $r\circ r=r$, i.e.
sc-smooth retractions. Whereas for Fr\'echet differentiability the image of such a retraction is easily shown to be a submanifold of $U$, the images of sc-smooth retractions can be much more general. Most strikingly, they can have locally varying dimensions. Of course, a good notion of differentiability comes with the chain rule so that from $r\circ r=r$ we deduce $Tr=(Tr)\circ (Tr)$. In  other words, the tangent map
of an sc-smooth retraction is again an sc-smooth retraction. If a subset $O$ of an sc-Banach space is the image of an sc-smooth retraction $r$, then $TO=Tr(TU)$  defines
the tangent space, and it turns out, that the definition does not depend on the choice of $r$. So we obtain quite general subsets of Banach spaces which have tangent spaces.
An  sc-smooth map $f:O\rightarrow O'$, where $O\subset E$ and $O'\subset F$ are sc-smooth retracts, is a map such that
$f\circ r:U\rightarrow F$ is sc-smooth, where $r$ is an sc-smooth retraction onto $O$. As it turns out the definition does not depend on the choice of $r$. Further,  one verifies that $Tf:=T(f\circ r)|TO$  defines  a map $TO\rightarrow TO'$ between tangent spaces and  that the definition also does not depend on the choice of $r$.

\begin{figure}[htbp]
\mbox{}\\[2ex]
\centerline{\relabelbox
\epsfxsize 4.4truein \epsfbox{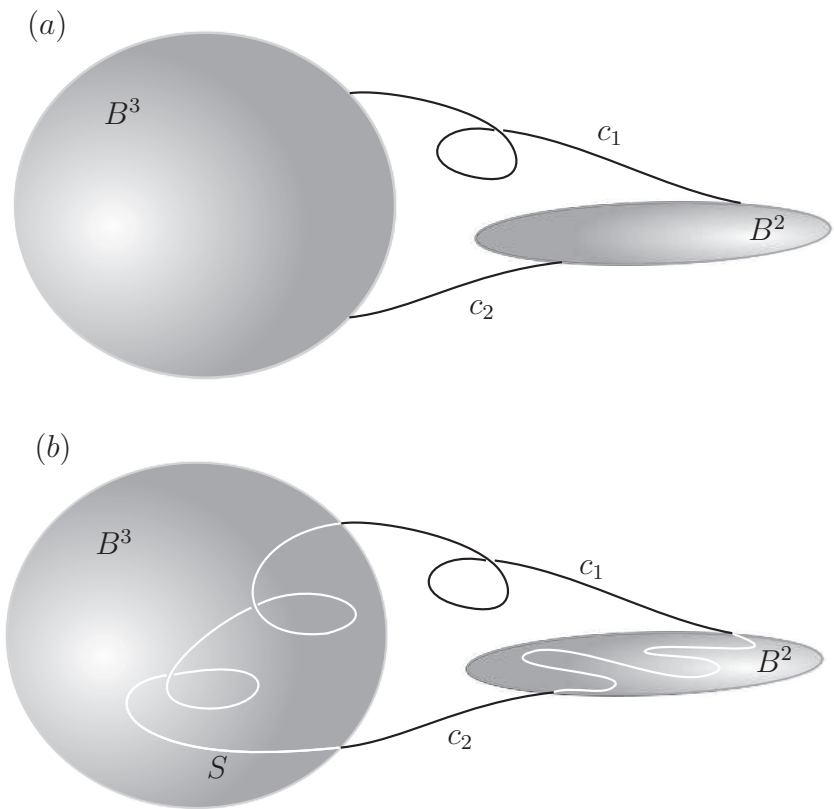}
\relabel {b1}{$B^3$}
\relabel {d1}{$B^2$}
\relabel {b}{$B^3$}
\relabel {d}{$B^2$}
\relabel {s}{$S$}
\relabel {c}{$c_1$}
\relabel {e2}{$c_2$}
\relabel {c1}{$c_1$}
\relabel {e1}{$c_2$}
\relabel {a}{$(a)$}
\relabel {a1}{$(b)$}
\endrelabelbox}
\caption{Figure a) shows a finite-dimensional M-polyfold $X$ which
is homeomorphic to the space consisting of the disjoint union of  an
open three-ball $B^3$ and an open two-ball $B^2$ connected by two
curves $c_1$, $c_2$. Figure b) shows the same M-polyfold containing
a one-dimensional $S^1$-like submanifold $S$. This submanifold could
arise as the zero set  of a transversal section of a strong
M-polyfold bundle $Y$ over $X$, which has varying dimensions.
Namely, over the three-ball  it is two-dimensional, over the two-disk
one-dimensional and otherwise it is trivial. The polyfold theory
would then guarantee a  natural smooth structure on the  solution
set $S$. }\label{porkbarrel}
\end{figure}

In summary, once we have a good notion of differentiability for maps between open sets, we also obtain a notion of differentiability for maps between smooth retracts. However, for the usual notion of differentiability, smooth retracts
are manifolds and one does not obtain anything beyond the usual differential geometry and its standard generalization to Banach manifolds. On the other hand sc-differentiability
opens up new possibilities with serious applications.
We generalize differential geometry by generalizing the notion of a manifold to that of an  M-polyfold. These  are metrizable spaces which are locally homeomorphic
to retracts with sc-smooth transition maps.  The theory as described in this paper even gives new objects in finite dimensions, see Figure \ref{porkbarrel}.

The generalization of an orbifold will be the notion of a polyfold. Then we define the notion of strong bundles over M-polyfolds and more generally polyfolds. The beautiful fact is that the compactified moduli spaces in Gromov-Witten, Floer-Theory and SFT can be thought of as zero-sets of sections of strong bundles
over polyfolds. These sections have suitable properties so that they are examples of the (generalized) Fredholm sections which will be described in this paper.

Assume that we have a strong bundle $W\rightarrow X$ over a polyfold with its sc-smooth structure.
We can define a particular class of sections $f$ which we shall refer to as Fredholm sections.  Moreover, one can introduce, as a consequence of
intrinsic properties of the spaces, the notion of an $\ssc^+$-multi-section. The latter is a map $\lambda:W\rightarrow {\mathbb Q}^+$ with suitable properties. Given $f$ and $\lambda$
we define the solution set of the pair $(f,\lambda)$ by
$$
{\mathcal S}=\{x\in X\ |\ \lambda(f(x))>0\}.
$$
The latter generalizes the following, more common, situation. Assume we consider instead of $f=0$ the perturbed problem
$f(x)=s(x)$. If we define $\lambda_s$ by $\lambda_s(w)=0$ for $w\not\in \hbox{graph}(s)$, and $\lambda(w)=1$ otherwise,
we can rewrite the problem $f(x)=s(x)$ as finding the solutions of $\lambda_s(f(x))>0$. Our approach generalizes this familiar situation.

We shall define what it means that $(f,\lambda)$ is in general position \textcolor{red}{and} in good position to the boundary. We shall also develop a perturbation theory
to bring solution sets into favorable positions.
There is also a notion of an  sc-smooth differential form on $X$, and for a generic solution set integration is defined. Even Stokes' theorem holds.
This can be used to define invariants. All this will be explained in more detail in the following.\\

\noindent{\bf Acknowledgement:}  The author thanks the Clay Institute for the opportunity to present this material at the 2008 Clay Research Conference.
Some of the work was done during a sabbatical at the Stanford mathematics department, and was supported in part by the American Institute of Mathematics. Thanks to P. Albers, Y. Eliashberg, E. Ionel, K. Wysocki and E. Zehnder for many stimulating discussions.

\section{Smoothness}
We begin with a question: Are there subsets of Banach spaces which are not submanifolds,
which however still carry something like a smooth structure and which can serve
as new local models for spaces? The answer is `yes', but one has to relax what one means by smoothness in infinite
dimensions. There is no choice, of course, what smoothness means in finite dimensions.
\subsection{Sc-Structures and Sc-Smooth Maps}
In interpolation theory, \cite{Tr}, given two Banach spaces $E\subset F$, general methods are developed to construct Banach spaces which interpolate
between $E$ and $F$. We take the concept of a scale (with suitable properties) from interpolation theory, but give it a new interpretation as a generalization of a smooth structure.
This is carried out in detail in \cite{HWZ1}.
\begin{defn}\label{d1}
Let $E$ be a Banach space. An {\bf sc-structure} for $E$ consists
of a nested sequence of Banach spaces $E_0\supset E_1\supset E_2\supset \cdots$ with $E_0=E$ so that
\begin{itemize}
\item[(1)]  The inclusion operator $E_{i+1}\rightarrow E_i$ is  a compact operator.
\item[(2)] $E_\infty=\bigcap_i E_i$ is dense in every $E_m$.
\end{itemize}
\end{defn}

As already said, such sc-structures occur in interpolation theory and are special cases of scales.  The interpretation
of an sc-structure as a generalization of a smooth structure, which we give soon, seems to be new.

\begin{examp} A typical example is $E=L^2({\mathbb R})$
with the sc-structure given by $E_m:=H^{m,\delta_m}({\mathbb R})$, where $H^{m,\delta_m}({\mathbb R})$ is the Sobolev space
of functions in $L^2$ so that the derivatives up to order $m$ multiplied by $e^{\delta_m |s|}$ belong to $L^2$.
Here $\delta_m$ is a strictly increasing sequence starting with $\delta_0=0$.
\end{examp}

If $E$ and $F$ are sc-Banach spaces, then $E\oplus F$ has a natural sc-structure given by
$$
(E\oplus F)_m=E_m\oplus F_m.
$$
Let us note that every finite-dimensional vector space has a unique sc-structure, namely the constant one,
where $E_i=E$. If $E$ is infinite-dimensional the constant sequence violates (1) of Definition \ref{d1}.

\begin{rem}
If $E$ is equipped with an sc-structure,
think of $x$ being a point in $E_i$, as a statement about the regularity of the point. A point in $E_\infty$ will be called a {\bf smooth point}.
In the classical theory the tangent space of $E$ at a point $x$ has a natural identification with $E$. This has to be modified in the sc-context.
Points $x\in E\setminus E_1$ do not have a tangent space, whereas points in $E_1$ have a tangent space which can be naturally identified with $E_0=E$.
Hence the tangent bundle $TE$ of $E$, which in the classical theory is $$
TE=E\oplus E,
$$
is in the sc-theory
$$
TE=E_1\oplus E_0=E_1\oplus E.
$$
 This space again has a natural
sc-structure given by $(TE)_m=E_{m+1}\oplus E_m$. Then we obtain the iterated tangent $T(TE)=E_2\oplus E_1\oplus E_1\oplus E_0$ and so on. We shall discuss this later in more detail, but note for the moment that the tangent bundle is defined
over $E_1$, i.e. $E_1\oplus E_0\rightarrow E_1$. Hence starting with the sc-Banach space $E$ its tangent bundle is $TE:=E_1\oplus E_0$ with the already
discussed sc-structure. In the previous example we therefore have
$$
TL^2({\mathbb R})=H^{1,\delta_1}({\mathbb R})\oplus L^2({\mathbb R})\ \ \hbox{and}\ \ (TL^2({\mathbb R}))_i=H^{i+1,\delta_{i+1}}({\mathbb R})\oplus H^{i,\delta_{i}}({\mathbb R}).
$$
\end{rem}
We continue with some considerations about linear sc-theory.
\begin{defn}
Let $E$ be an sc-Banach space and $F\subset E$ a linear subspace. We call $F$ an {\bf sc-subspace} provided the filtration $F_i=F\cap E_i$
turns $F$ into an sc-Banach space. If $F\subset E$ is an sc-Banach space, then we say that it has an sc-complement, provided there exists an
sc-subspace $G$ such $F_i\oplus G_i=E_i$ as topological linear sum for all $i$.
\end{defn}

Let us note that a finite-dimensional subspace $F$ of $E$ has an sc-complement if and only if $F\subset E_\infty$, see \cite{HWZ1}.

The linear operators of interest are those linear operators $T:E\rightarrow F$,
which map $E_m$ into $F_m$ for all $m$, such that $T:E_m\rightarrow F_m$ is a bounded linear operator.
We call $T$ an {\bf sc-operator}. An {\bf sc-isomorphism} $T:E\rightarrow F$ is a bijective sc-operator so that its inverse is also an sc-operator.
Of particular interest are the linear {\bf  sc-Fredholm operators}:

\begin{defn}
A sc-operator $T:E\rightarrow F$ is said to be {\bf sc-Fredholm} provided there exist sc-splittings $E=K\oplus X$ and $F=Y\oplus C$ so that $C$ is finite-dimensional, $Y=T(X)$, and $T:X\rightarrow Y$ defines a linear sc-isomorphism.
\end{defn}
We note that the above implies that $E_m=X_m\oplus K$ and $F_m=T(X_m)\oplus C$ for all $m$.  The Fredholm index is by definition
$$
\ind(T)=\dim(K)-\dim(C).
$$
Let us also observe that for every $m$ we have a linear Fredholm operator (in the classical sense)  $T:E_m\rightarrow F_m$, which in particular have the same index and identical kernels.

Next we begin with the preparations to introduce the notion of an sc-smooth map.
\begin{defn}
A {\bf partial quadrant} $C$ in an sc-Banach space $E$ is a closed convex subset with the property
that there exists an sc-Banach space $W$ and an sc-isomorphism $T:E\rightarrow {\mathbb R}^n\oplus W$
mapping $C$ onto $[0,\infty)^n\oplus W$.
\end{defn}
The main point for introducing partial quadrants is the needed generality of our theory, which in its applications to SFT has to provide spaces with boundary with corners.

We note that $C=E$ is a partial quadrant in $E$. In the following we are interested in relatively open subsets
of partial quadrants  $C$.
If $U\subset C\subset  E$ is a relatively open subset we can define an sc-structure for $U$ by the nested sequence
$U_i$ given by $U_i=E_i\cap U$. We note that $U_\infty=\bigcap U_i$ is dense in every $U_m$.
Considering $U$ with its sc-structure we see that $U_{i_0}$ also admits an sc-structure defined by
$$
{(U_{i_0})}_m:=U_{i_0+m}.
$$
We write $U^{i_0}$ for $U_{i_0}$ equipped with this sc-structure.
Given two such sc-spaces $U$ and $V$ we write $U\oplus V$ for $U\times V$ equipped with the
obvious sc-structure. Now we can give the rigorous definition of the tangent $TU$ of a relatively open subset $U$ of a partial quadrant $C$ in an sc-Banach space $E$.

\begin{defn}
The {\bf tangent} $TU$ of a relatively open subset $U\subset C\subset E$ of the sc-Banach space $E$ is defined by $TU=U^1\oplus E$.
\end{defn}
We note that
$$
(TU)_i=U_{1+i}\oplus E_i.
$$
Recalling our previous example we have
$$
TL^2({\mathbb R}) = H^{1,\delta_1}({\mathbb R})\oplus L^2({\mathbb R})\ \ \hbox{and}\ \
(TL^2({\mathbb R}))_i = H^{i+1,\delta_{i+1}}({\mathbb R})\oplus H^{i,\delta_i}({\mathbb R}).
$$

Given two relatively open subsets $U$ and $V$ of partial quadrants, a map $f:U\rightarrow V$ is said to be of class ${\bf sc}^\mathbf{0}$ provided
for every $m$ the map $f$ maps $U_m$ into $V_m$ and the map $f:U_m\rightarrow V_m$ is continuous. The following example takes a little bit of work.
\begin{examp}
Take $L^2({\mathbb R})$ with the previously defined sc-structure and define
$$
\Phi:{\mathbb R}\oplus L^2({\mathbb R}) \rightarrow L^2({\mathbb R}):(t,u)\rightarrow \Phi(t,u),
$$
where $\Phi(t,u)(s)=u(s+t)$. Then $\Phi$ is $\ssc^0$. As we shall see later, this will also be an example for an sc-smooth map.
\end{examp}

Next we define the notion of an $\ssc^1$-map.

\begin{defn}
Let $U\subset C\subset  E$ and $V\subset D\subset F$ be relatively open subsets of partial quadrants in sc-Banach spaces.
An  $\ssc^0$-map $f:U\rightarrow V$ is said to be ${\bf sc}^\mathbf{1}$ provided for every $x\in U_1$ there exists a continuous linear operator $Df(x):E_0\rightarrow F_0$ so that the following holds.
\begin{itemize}
\item[(1)] For $h\in E_1$ with $x+h\in U$ we have
$$
\lim_{\parallel h\parallel_1\rightarrow 0} \frac{1}{\parallel h \parallel_1}\cdot \parallel f(x+h)-f(x)-Df(x)h\parallel_0 =0.
$$
\item[(2)] The map $Tf$ defined by $Tf(x,h)=(f(x),Df(x)h)$ for $(x,h)\in TU$ defines an  $\ssc^0$-map $Tf:TU\rightarrow TV$.
    \end{itemize}
    \end{defn}

Inductively we can define what an  $sc^k$ map is, and what it means to be $sc^\infty$. The following result shows that the chain-rule holds.
\begin{thm}[Chain-Rule]
Assume that $U,V$ and $W$ are relatively open subsets of partial quadrants and $f:U\rightarrow V$ and $g:V\rightarrow W$
are $\ssc^1$-maps. Then $g\circ f$ is $\ssc^1$ and $T(g\circ f)=(Tg)\circ (Tf)$.
\end{thm}
\begin{examp}
One can show, see \cite{H2} (taken from \cite{HWZ-polyfolds1}), that
the map $\Phi$ from the previous example is sc-smooth.
\end{examp}

\subsection{Sc-Smooth
 Spaces and M-Polyfolds}

Now we are in the position to introduce new local models for smooth spaces.
The following presentation is taken from \cite{H2}\footnote{This lecture note and \cite{H3} are based on graduate courses given at Courant in 2004/2005 and Stanford in 2008.} and generalizes somewhat
the results in \cite{HWZ1}\footnote{These generalizations grew out of discussions with my graduate students and postdocs, most notably Joel Fish and Peter Albers, while running a seminar at Courant on \cite{HWZ1,HWZ2,HWZ3,HWZ5}. All of these generalizations follow immediately from the proofs
in \cite{HWZ1,HWZ2,HWZ3}. The benefit of the current version is that the formalism gets somewhat easier and more transparent.}.
The interesting thing about sc-smoothness is the fact that there are many smooth retractions
with complicated images, so that one obtains a large 'library' of local models for spaces, which is sufficient to describe problems
occurring when studying partial differential equations with analytical limiting behavior allowing for bubbling-off and similar analytical phenomena.

\begin{defn}
Let $U\subset C\subset E$ be a relatively open subset in a partial quadrant of the sc-Banach space $E$.
A map $r:U\rightarrow U$ is called a {\bf $\ssc^\infty$-retraction} provided it sc-smooth and $r\circ r=r$.
\end{defn}
The chain rule implies that for a $\ssc^\infty$-retraction $r$ its tangent map $Tr$ is again an $\ssc^\infty$-retraction.
We call the image $O=r(U)$ of an $\ssc^\infty$-retraction $r:U\rightarrow U$ an {\bf $\ssc^\infty$}-{\bf retract}.

The crucial definition is the following.
\begin{defn}
A {\bf local sc-model} is a triple $(O,C,E)$, where $E$ is an sc-Banach space, $C$ a partial quadrant and $O\subset C$ an  $\ssc^\infty$-retract
given as the image of an  sc-smooth retraction $r:U\rightarrow U$ defined on a relatively open subset  $U$ of $C$.
\end{defn}

The following lemma is easily established, see \cite{H2}.
\begin{lem}
Assume that $(O,C,E)$ is a local sc-model and $r$ and $s$ are sc-smooth retractions defined on relatively open subsets $U$ and $V$
of $C$, respectively,  having $O$ as the image. Then $Tr(TU)=Ts(TV)$.
\end{lem}

In view of this lemma we can define the tangent of a local sc-model which again is a local sc-model as follows.
\begin{defn}
The {\bf tangent} of the local sc-model $(O,C,E)$ is defined by
$$
T(O,C,E):= (TO,TC,TE),
$$
where $TO=Tr(TU)$ for any sc-smooth retraction $r:U\rightarrow U$ having $O$ as the image,
where $U$ is relatively open in $C$.
\end{defn}
\begin{rem}Let us observe that if $(O,C,E)$ is a local sc-model and $O'$ an open subset of $O$, then $(O',C,E)$ is again a local sc-model.
Indeed, if $r:U\rightarrow U$ is an  sc-smooth retraction, where $U$ is relatively open in $C$, with $O=r(U)$, then define $U'=r^{-1}(O')$.
This is relatively open in $C$ and $r'=r|U':U'\rightarrow U'$ is an sc-smooth retraction with image $O'$.
\end{rem}

A map $f:O\rightarrow O'$ between two local sc-models is {\bf sc-smooth (or $\ssc^k$)} provided $f\circ r:U\rightarrow E'$ is sc-smooth (or $\ssc^k$). One easily verifies that the definition
does not depend on the choice of $r$. We can define the tangent map $Tf:TO\rightarrow TO'$ by
$$
Tf:= T(f\circ r)|Tr(U).
$$
As it turns out this is well-defined and does not depend on the choice of $r$ as long as it is compatible with $(O,C,E)$.

\begin{thm}[Chain Rule]
Assume that $(O,C,E)$, $(O',C',E')$ and $(O'',C'',E'')$ are local sc-models and $f:O\rightarrow O'$ and $g:O'\rightarrow O''$
are $\ssc^1$. Then $g\circ f:O\rightarrow O''$ is $\ssc^1$ and
$$
T(g\circ f) =(Tg)\circ (Tf).
$$
\end{thm}
The following Remark \ref{2.18} explains how the current account is related to \cite{HWZ1,HWZ2,HWZ3,HWZ4}.
\begin{rem}\label{2.18} In the series of papers \cite{HWZ1,HWZ2,HWZ3,HWZ4} we developed a generalized Fredholm theory in a slightly more restricted situation, which is more than enough for the applications.  Namely rather than considering $\ssc$-smooth retractions and $\ssc$-smooth retracts, splicings and open subsets of splicing cores were considered, which one can view as a special case.
Namely a splicing consists of a relatively open subset $V$ of a partial quadrant $C$ in some sc-Banach space $W$ and a family of bounded linear projections
$\pi_v:E\rightarrow E$, $v\in V$, where $E$ is another sc-Banach space, so that the map
$$
V\oplus E\rightarrow E:(v,e)\rightarrow\pi_v(e)
$$
is sc-smooth. Then the associated splicing core is $K$, defined by
$K=\{(v,e)\in V\oplus E \ |\ \pi_v(e)=e\}$. Clearly $V\oplus E$ is a relatively open subset of the partial quadrant $C\oplus E$ in $W\oplus E$ and $r(v,e):=(v,\pi_v(e))$ is an sc-smooth retraction. The associated retract is, of course, the splicing core $K$. If $O$ is an open subset of $K$ we know that it is again an sc-smooth retract.  Let us note that in all our applications the retractions are obtained from splicings.
\end{rem}
We demonstrate first how the definition of a manifold can be generalized.
Let $Z$ be a metrizable topological space. A {\bf chart} for $Z$ is a tuple $(\varphi,U,(O,C,E))$, where $\varphi:U\rightarrow O$ is a homeomorphism
and $(O,C,E)$ is a local sc-model. We say that two such charts are sc-smoothly compatible provided
$$
\psi\circ\varphi^{-1}:\varphi(U\cap V)\rightarrow \psi(U\cap V)
$$
is sc-smooth and similarly for $\varphi\circ \psi^{-1}$. Here $(\psi,V,(P,D,F))$ is the second chart. Note that the sets $\varphi(U\cap V)$ and $\psi(U\cap V)$
are sc-smooth retracts for sc-smooth retractions defined on relatively open sets in $C$ and $D$, respectively. An sc-smooth atlas for  $Z$
consists of a family of sc-smoothly compatible charts so that their domains cover $Z$. Two sc-smooth atlases are compatible provided their union is an sc-smooth atlas. This defines
an equivalence relation.

\begin{defn}
Let $Z$ be a metrizable space. An  {\bf sc-smooth structure} on $Z$ is given by an sc-smooth atlas. Two sc-smooth structures
are {\bf equivalent} if the union of the two associated atlases defines again an sc-smooth structure.
An  {\bf sc-smooth space} is a metrizable space $Z$ together with an equivalence class of sc-smooth structures.
\end{defn}
We note that these spaces have a natural filtration $Z_0\supset Z_1\supset Z_2\supset \cdots $. The points in $Z_i$ one should view as the points
of some regularity $i$.
The sc-smooth spaces are a very general type of space on which one can define sc-smooth functions.

 It is possible to generalize many of the constructions
from differential geometry to these spaces. If we have an sc-smooth partition of unity we can define Riemannian metrics and consequently a curvature tensor.
Note however that curvature would only be defined  at points of regularity at least $2$. The existence of an sc-smooth partition of unity depends on the sc-structure\footnote{Recall that the existence of a smooth partition of unity on a Banach space is a classical question and related to differentiability questions of norms.}.

The {\bf tangent space} at a point of level at least one is defined in the same way as one defines them for Banach manifolds, see \cite{LA}. Namely one considers
tuples $(z,\varphi,U,(O,C,E),h)$, where $z\in Z_1$ and $(\varphi,U,(O,C,E))$ a chart, so that $z\in U$, and $h\in T_{\varphi(z)}O$. Two such tuples, say the second is $(z',\varphi',U',(O',C',E'),h')$, are said to be equivalent
provided $z=z'$ and $T(\varphi'\circ \varphi^{-1})(\varphi(z))h=h'$.
An equivalence class $[(z,\varphi,U,(O,C,E),h)]$, then is a tangent vector at $z$. The tangent space at $z\in Z_1$ is denoted by $T_zZ$ and we define
$TZ$ as
$$
TZ=\bigcup_{z\in Z_1} \{z\}\times T_zZ.
$$
One can show that $TZ$ has a natural sc-smooth structure so that the natural map $TZ\rightarrow Z^1$ is sc-smooth.

For applications there are two particular examples of sc-smooth spaces which are very useful.
\begin{defn}
An  sc-smooth space, which admits an atlas where the charts have the form $(\varphi,U,(O,E,E))$, is called a {\bf M-polyfold} (without boundary).
\end{defn}

\begin{examp} Consider the metrizable space $Z$ given as the subspace of ${\mathbb R}^2$ defined by
$$
Z=\{(s,t)\in {\mathbb R}^2\ |\ t=0\ \hbox{if}\ s\leq 0\}.
$$
Then $Z$ admits the structure of an M-polyfold without boundary. In order to see this, one constructs a topological embedding into ${\mathbb R}\oplus L^2({\mathbb R})$, where
$L^2({\mathbb R})$ has the previously introduced sc-structure, in such a way that the image is an sc-smooth retract. Here the idea of an sc-smooth splicing comes in handy! Take a smooth, compactly supported map $\beta:{\mathbb R}\rightarrow [0,\infty)$ with $\int \beta(t)^2 ds=1$.
Denote by $f_s$, for $s\in (0,\infty)$ the unit length element in $L^2$ defined by
$$
f_s(t)=\beta(t+e^\frac{1}{s}).
$$
For $s\in (-\infty,0]$ we define $f_s=0$. Then let $\pi_s$ be the $L^2$-orthogonal projection onto the subspace spanned by $f_s$.
Then a somewhat lengthy computation, carried out in \cite{HWZ-polyfolds1}, shows that
$$
r:{\mathbb R}\oplus L^2\rightarrow {\mathbb R}\oplus L^2:r(s,u)=(s,\pi_s(u))
$$
is an $\ssc$-smooth retraction, with obvious image $O$ being
$$
\{(s,t\cdot f_s)\ |\ (s,t)\in {\mathbb R}^2\}.
$$
Hence $(O,{\mathbb R}\oplus L^2,{\mathbb R}\oplus L^2)$ is a local sc-model.
We note that it has varying dimension.
The map
$$
Z\rightarrow {\mathbb R}\oplus L^2:(s,t)\rightarrow (s, t\cdot f_s)
$$
is a homeomorphic embedding onto $O$. The map is clearly continuous and
injective and has image $O$.  Define ${\mathbb R}\oplus L^2\rightarrow {\mathbb R}^2$ by
$$
(s,x)\rightarrow (s,\int_{\mathbb R} x(t)f_s(t)dt).
$$
This map is continuous and its restriction to $O$ is the inverse of the previously defined map.
Hence we obtain the structure of an  M-polyfold on $Z$. This gives us the first example of a finite-dimensional space, with varying dimension, which has
a generalized manifold structure. We also note that the induced filtration is constant, so that a tangent space is defined at all points. This is due to the fact that the local model $O$ lies entirely in the smooth part of ${\mathbb R}\oplus L^2$. It is instructive to study sc-smooth curves
$\phi:(-\varepsilon,\varepsilon)\rightarrow O$ satisfying $\phi(0)=(0,0)$.
\end{examp}

We also would like to generalize the notion of a manifold with boundary with corners. This forces us to allow  only a certain class
of local sc-models. Let us observe that in the definition of a local sc-model $(O,C,E)$ we did not require $O$ to be in a particular position to the boundary of $C$. To define the boundary of $O$ as $\partial O:= O\cap \partial C$ is for general local sc-models not(!) a good idea!

However, things change, if we take triples $(O,C,E)$ with more structure.
For a point $x\in C\subset E$ we denote by $d_C(x)$ the {\bf degeneracy index} which is defined as follows. Take an sc-isomorphism $T:E\rightarrow {\mathbb R}^n\oplus W$
mapping $C$ to $[0,\infty)^n\oplus W$. Then $d_C(x)$ is the number of indices $i\in \{1,\ldots , n\}$ such that $r_i=0$, where $T(x)=(r_1,\ldots ,r_n,w)$.

\begin{defn}
Let $U\subset C\subset E$ be a relatively open subset of the partial quadrant $C$. An  sc-smooth retraction $r:U\rightarrow U$ is said to be {\bf neat} provided
the following holds:
\begin{itemize}
\item[(1)] For every smooth point $x\in O$ the kernel $N$ of $Id-Dr(x)$ has an  sc-complement $M$ contained in $C$.
\item[(2)] For every point $x\in O$ there exists a sequence of smooth points $(x_k)$ in $O$ converging to $x$
and satisfying $d_C(x_k)=d_C(x)$.
\end{itemize}
\end{defn}

Now we can define the local sc-models of interest.
\begin{defn}
A local sc-model $(O,C,E)$ is called {\bf neat} provided there exists a neat sc-smooth retraction $r:U\rightarrow U$,
where $U$ is relatively open in $C$ with $O=r(U)$.
\end{defn}

We define the {\bf degeneracy index} $d_O$ for a neat sc-model by $d_O:O\rightarrow {\mathbb N}$ and $d_O=d_C|O$.
The following crucial theorem tells us that the degeneracy index is an invariant under sc-diffeomorphism.
\begin{thm}[Corner Recognition] Assume we are given neat local sc-models $(O,C,E)$ and $(O',C',E')$  and an sc-diffeomorphism $f:O\rightarrow O'$. Then $d_O(x)=d_{O'}(f(x))$ for all $x\in O$.
\end{thm}
The proof of this theorem is given in \cite{H2} and is an easy adaption of
a similar result in \cite{HWZ1} dealing with a splicing-based differential geometry.
\begin{defn}
An  sc-smooth space build on neat local sc-models is called a {\bf M-polyfold} with boundary with corners.
\end{defn}
If $Z$ is an  M-polyfold with boundary with corners, then we have a well-defined degeneracy map $d:Z\rightarrow {\mathbb N}$.
The points $x$ with $d(x)=0$ are the interior points, whereas the points with $d(x)\geq 1$ are the boundary points. If $d(x)\geq 2$
we have a corner of order $d(x)$.

For the later discussions about Fredholm theory, we need the notion of a submanifold $M$ of an  M-polyfold $X$. We need the following definition
\begin{defn}\label{good-position}
Let $C$ be a partial quadrant in the sc-Banach space $E$.
The closed linear subspace $N$ of $E$ is in {\bf good position} to the partial quadrant $C$ if $N\cap C$ has a nonempty interior in $N$ and there exists an sc-complement denoted by $N^\perp$ in $E$ and a constant $c>0$ so that for every vector $(n,m)\in N\oplus N^\perp$ satisfying
    $$
    \norm{m}_E\leq c\cdot\norm{n}_E
    $$
    the statements $n+m\in C$ and $ n\in C$ are equivalent. If this equivalence holds for an sc-complement $N^\perp$, we call it a {\bf good complement}.
\end{defn}
A submanifold of an  M-polyfold can now be defined.
\begin{defn}
Let $X$ be an M-polyfold and $M$ a subset equipped with the induced topology. The subset $M$ is called a {\bf finite dimensional submanifold} of $X$ provided the following holds.
\begin{itemize}
\item[(1)] The subset $M$ lies in $X_\infty$.
\item[(2)] At every point $m\in M$, there exists an M-polyfold chart
$(\varphi,U,(O,C,E))$ with $m\in U$ so that $\varphi(m)=0$ having the following property. Putting
$M^\ast=\varphi(M\cap U)$, there exists a finite-dimensional smooth linear subspace $N\subset E$ in good position to $C$, a corresponding good sc-complement $N^\perp$,  an open neighborhood $Q$ of $0\in C\cap N$, and an sc-smooth map $A:Q\rightarrow N^\perp$, with $DA(0)=0$ and $A(0)=0$, so that the map
$$
\Gamma:Q\rightarrow E:q\rightarrow q+A(q)
$$
has as image precisely $M^\ast$.
\item[(3)] The map $\Gamma:Q\rightarrow M^\ast$ is a homeomorphism
\end{itemize}
The map $\Phi:=\varphi^{-1}\circ\Gamma$ is called a good parametrization of a neighborhood of $m$.
\end{defn}
It has been proved in \cite{HWZ2}, that any two good parameterizations are smooth\-ly compatible. Hence a finite-dimensional submanifold  $M$ of an M-polyfold carries in a natural way the structure of a smooth manifold with boundary with corners ($\partial M=\emptyset$ is, of course, possible).

In a later section we shall introduce polyfolds, which are more general than M-polyfold with boundary with corners. These will be the spaces needed for symplectic field theory.

The following example is relevant for Morse-homology and can be modified
and extended to Floer-theory.
\begin{examp}
Let $a,b$ and $c$ be three (mutually) distinct points in ${\mathbb R}^n$.
Fix two smooth maps $\varphi:{\mathbb R}\rightarrow {\mathbb R}^n$
and $\psi:{\mathbb R}\rightarrow {\mathbb R}^n$ such that
$$
\varphi(s)=a,\ \psi(s)=b\ \hbox{for}\ \ s<<0
$$
and
$$
\varphi(s)=b,\ \psi(s)=c\ \hbox{for}\ \ s>>0.
$$
Finally we take a third smooth function $\sigma$, which in a similar way connects $a$ with $c$.
If $H^2=H^2({\mathbb R},{\mathbb R}^n)$ denotes the usual Sobolev space
we define $X'(a,b)=\varphi+H^2$ and $X'(b,c)=\psi+H^2$ and $X'(a,c)=\sigma+H^2$. All these spaces have natural topologies and also admit natural ${\mathbb R}$-actions by the action of ${\mathbb R}$ on itself. Also the spaces do not depend on the particular choices of the smooth functions $\varphi,\psi$ and $\sigma$, provided they have the required properties. Denote by $X(a,b)$, $X(b,c)$ and $X(a,c)$ the corresponding quotient spaces. We call an element of the latter spaces a `path'. These spaces are second countable and paracompact, in fact metrizable. Define the set $\bar{X}$ by
$$
\bar{X}= X(a,c)\coprod (X(a,b)\times X(b,c)).
$$
It is the space of possibly once-broken paths from $a$ to $c$.
One can define on $\bar{X}$ a metrizable topology in a natural way inducing on
$X(a,c)$ and $X(a,b)\times X(b,c)$ the already given topology, so that in  addition $\bar{X}$ is connected and $X(a,c)$ is open and dense. This topology captures in a geometric way how a path from $a$ to $c$ decomposes
to a broken path and defines what it means that an unbroken path is close to a broken path. Assume we are given a sequence $0=\delta_0<\delta_1<\delta_2<\cdots$ of increasing weights. Denote by $\bar{X}_m$ the subset of $\bar{X}$ obtained by replacing $H^2$ by the weighted Sobolev space $H^{2+m,\delta_m}$, where $x$ belongs to this space provided all derivatives up to order $m+2$ weighted by $e^{\delta_m |s|}$ belong to $L^2$. One can show that $\bar{X}$ has in a natural way the structure
of an  M-polyfold with boundary with corners for which the level $m$ subset of $\bar{X}$ corresponds
to $\bar{X}_m$. We have a well-defined degeneracy map
$$
d:\bar{X}\rightarrow \{0,1,2,\ldots \}
$$
which only takes the values $0$ and $1$. Namely $d|X(a,c)\equiv 0$ and
$d|(X(a,b)\times X(b,c))\equiv 1$, so that the broken paths are indeed
the boundary points.
Clearly, a similar statement for maps defined on cylinders approaching loops in a symplectic manifold, makes the example relevant for Floer-theory, see \cite{H2} and\cite{HWZ-polyfolds1} for details.
\end{examp}

\subsection{Strong Bundles}
The notion of a strong bundle is designed to give additional structures in the Fredholm theory, which guarantee a compact perturbation and transversality theory.
The crucial point is the fact that there will be a well-defined vector space of perturbations, which have  certain compactness properties, but on the other
hand are plentiful enough to allow for all kinds of perturbations needed for different versions of the Sard-Smale theorem, \cite{smale}, in the Fredholm theory.

Let us start with a non-symmetric product $U\triangleleft F$, where $U$ is a relatively open subset in some partial quadrant of the sc-Banach space $E$, and $F$ is also an sc-Banach space.
By definition, as a set $U\triangleleft F=U\times F$, but in addition it has a double filtration
$$
(U\triangleleft F)_{m,k}=U_m\oplus F_k
$$
defined for all pairs $(m,k)$ satisfying $0\leq k\leq m+1$. We view $U\triangleleft F\rightarrow U$ as a bundle with base space $U$ and fiber $F$,
where the double filtration has the interpretation that above a point $x\in U$ of regularity $m$ it makes sense to talk about fiber regularity
of a point $(x,h)$ up to order $k$ provided $k\leq m+1$.
We shall explain this asymmetry and range of parameters in the filtration below. Given $U\triangleleft F$, we might consider the associated sc-spaces $U\oplus F$ and $U\oplus F^1$.

Of interest  for  us are the maps
$$
\Phi:U\triangleleft F\rightarrow V\triangleleft G
$$
of the form
$$
\Phi(u,h)=(\varphi(u),\phi(u,h))
$$
which are linear in $h$. We say that the map is of class $\ssc^0_\triangleleft$ provided it induces $\ssc^0$-maps
$U\oplus F^i\rightarrow V\oplus G^i$ for $i=0,1$.
We define the tangent $T(U\triangleleft F)$ by
$$
T(U\triangleleft F) = (TU)\triangleleft (TF).
$$
Note that the order of the factors is different from the order in $T(U\oplus F)$. One has to keep this in mind. Indeed,
$$
T(U\triangleleft F)=U_1\oplus E\oplus F_1\oplus F\ \ \hbox{and}\ \ T(U\oplus F)=U_1\oplus F_1\oplus E\oplus F.
$$
A map $\Phi:U\triangleleft F\rightarrow V\triangleleft G$ is of class $\ssc^1_\triangleleft$ provided
the maps $\Phi:U\oplus F^i\rightarrow V\oplus G^i$ for $i=0,1$ are $\ssc^1$. Taking the tangents of the latter,
gives after rearrangement, the $\ssc^0_\triangleleft$-map
$$
T\Phi:(TU)\triangleleft (TF)\rightarrow (TV)\triangleleft (TG).
$$
Iteratively we can define what it means that a map is $\ssc^k_\triangleleft$ for $k=1,2,\ldots$ and we can also define $\ssc_\triangleleft$-smooth maps.

Given $U\triangleleft F\rightarrow U$, an sc-smooth section $f$ is map of the form
$x\rightarrow (x,\bar{f}(x))$ such that the induced map $U\rightarrow U\oplus F$ is sc-smooth. In particular, $f$ is `horizontal' with respect to the filtration, i.e. a point on level $m$ is mapped to a point of bi-level $(m,m)$. This can be considered as a convention, and it is precisely this convention which is responsible for the filtration constraint $k\leq m+1$. There is another class of sections called $\ssc^+$-sections. These are sc-smooth sections
of $U\triangleleft F\rightarrow U$ which induce  sc-smooth maps
$U\rightarrow U\oplus F^1$.  In particular, if $s$ is an $\ssc^+$-section of  $U\triangleleft F\rightarrow U$ and  $s(x)=(x,\bar{s}(x))$  for $x\in U_m$  then  $\bar{s}(x)\in F_{m+1}$.
This type of sections will be important for the perturbation theory.
Indeed, it is a kind of compact perturbation theory since the inclusion $F_{m+1}\rightarrow F_m$ is compact.
We give an example before we generalize an earlier discussion about retracts and retractions to bundles of the type  $U\triangleleft F\rightarrow U$.

\begin{examp}
Let us denote by $E$ the Sobolev space $H^1(S^1,{\mathbb R}^n)$ of loops.
We define an sc-structure by $E_m= H^{1+m}(S^1,{\mathbb R}^n)$. Further we define
$F=L^2(S^1,{\mathbb R}^n)=H^0(S^1,{\mathbb R}^n)$ which we filter via
$F_m=H^m(S^1,{\mathbb R}^n)$. Finally we introduce $E\triangleleft F\rightarrow E$. Then we can view the map $f:x\rightarrow\dot{x}$ as an sc-smooth section. In particular,  $f$ maps $E_m$ into $E_m\oplus F_m$.
We observe that the filtration of $F$ is picked in such a way that the first order differential operator $x\rightarrow \dot{x}$ is an sc-smooth section,
in particular,  it is horizontal, i.e. the choices are made in such a way that they comply with our convention that sc-smooth sections are index preserving.
The map $x\rightarrow x$ can be viewed as an $\ssc^+$-section. Then $x\rightarrow \dot{x}+x$ is an sc-smooth section obtained from the sc-smooth section $x\rightarrow \dot{x}$ via the perturbation by an $\ssc^+$-section.
Consider now a smooth vector bundle map
$$
\Phi:{\mathbb R}^n\oplus {\mathbb R}^n\rightarrow
{\mathbb R}^n\oplus {\mathbb R}^n
$$
of the form
$$
\Phi(x,h)=(\varphi(x),\phi(x)h),
$$
where $\varphi:{\mathbb R}^n\rightarrow {\mathbb R}^n$ is a diffeomorphism and for every $x\in {\mathbb R}^n$ the map $\phi(x):{\mathbb R}^n\rightarrow {\mathbb R}^n$ is a linear isomorphism. Then we define for $(x,h)\in E\oplus F$ the element $\Phi_\ast(x,h)(t)=(\varphi(x(t)),\phi(x(t))h(t))$.
Note that if $x\in E_m$ and $h\in F_{k}$ for $k\leq m+1$, then $\Phi_\ast(x,h)=:(y,\ell)$ satisfies $y\in E_m$ and $\ell\in F_k$.
However, if $x\in E_m$ and $y\in F_k$ for some $k>m+1$ we cannot conclude that $\ell\in F_k$. We can only say that $\ell\in F_{m+1}$. Now one easily verifies that
$$
\Phi_\ast:E\triangleleft F\rightarrow E\triangleleft F
$$
is $\ssc_\triangleleft$-smooth. This justifies our constraint $k\leq m+1$ for the double filtration. Even more is true, which however is irrelevant for us, namely $\Phi_\ast:E_m\oplus F_{m+i}\rightarrow E_m\oplus F_{m+i}$ is smooth in the classical sense for all $m$ and $i=0,1$.
\end{examp}

\begin{defn}
An  {\bf $\ssc^\infty_\triangleleft$-retraction} is an  $\ssc_\triangleleft$-smooth map
$$
R:U\triangleleft F\rightarrow U\triangleleft F
$$
with the property $R\circ R=R$.
\end{defn}
Of course, $R$ has the form $R(u,h)=(r(u),\phi(u,h))$ with $r$ being an $\ssc$-smooth retraction and $\phi(u,h)$ linear in the fiber.
Given $R$, we can define its image $K=R(U\triangleleft F)$ and $O=r(U)$. Then we have a natural map
$$
p:K\rightarrow O
$$
We may view this as the local model for a strong bundle. Observe that $K$ has a double filtration and $p$ maps points of regularity
$(m,k)$ to  points of regularity $m$.
\begin{defn}
The tuple $(K,C\triangleleft F,E\triangleleft F)$, where $K$ is a subset of $C\triangleleft F$, so that there exists an
$\ssc^\infty_\triangleleft$-retraction $R$ defined on $U\triangleleft F$, where $U$ is relatively open subset of $C$ and $K=R(U\triangleleft F)$,
is called a {\bf local strong bundle model}.
\end{defn}
Starting with $(K,C\triangleleft F,E\triangleleft F)$ we have the projection $K\rightarrow E$ and denote its image by $O$
and the induced map by $p:K\rightarrow O$. One can define $T(K,C\triangleleft F,E\triangleleft F)$ by
$$
T(K,C\triangleleft F,E\triangleleft F)=(TK, TC\triangleleft TF, TE\triangleleft TF),
$$
where $TK$ is the image of $TR$. As before we can show that the definition does not depend on the choice of $R$.

Now we are in the position to define the notion of a strong bundle. Let $p:W\rightarrow X$ be a surjective continuous map between two metrizable spaces,
so that for every $x\in X$ the space $W_x:=p^{-1}(x)$ comes with the structure of a Banach space.  A {\bf strong bundle chart} is a tuple
$(\Phi,p^{-1}(U),(K,C\triangleleft F,E\triangleleft F))$, where $\Phi:p^{-1}(U)\rightarrow K$ is a homeomorphism, covering a homeomorphism $\varphi:U\rightarrow O$,
which between each fiber is a bounded linear operator.  We call two such charts {\bf $\ssc_\triangleleft$-smoothly equivalent} if the associated transition maps
are $\ssc_\triangleleft$-smooth.  We can define the notion of a {\bf strong bundle atlas} and can define the notion of equivalence of two such atlases.
\begin{defn}
Let $p:W\rightarrow X$ be as described before. A {\bf strong bundle structure} for $p$ is given by a strong bundle atlas. Two strong bundle structures are {\bf equivalent}
if the associated atlases are equivalent. Finally $p$ equipped with an equivalence class of strong bundle atlases is called a {\bf strong bundle}.
\end{defn}

Let us observe that a strong bundle $p:W\rightarrow X$ admits a double filtration $W_{m,k}$ with $0\leq k\leq m+1$. By forgetting part of this double filtration we observe that $W(0)$,
which is $W$ filtered by $W(0)_m:=W_{m,m}$ has in a natural way the structure of an  M-polyfold. The same is true for $W(1)$ which is the space $W_{0,1}$ equipped with the filtration $W(1)_m:=W_{m,m+1}$. Obviously the maps $p:W(i)\rightarrow X$ for $i=0,1$ are sc-smooth.

The previously introduced notions of $\ssc$-smooth sections and $\ssc^+$-sections for $U\triangleleft F\rightarrow U$ generalize as follows.

A {\bf sc-smooth section} of a strong bundle is an  $\ssc^0$-map $s:X\rightarrow W$ with $p\circ s=Id_X$ such that $s:X\rightarrow W(0)$ is sc-smooth.
The vector space of all such sections is written as $\Gamma(p)$. Very important for the perturbation theory are the so-called $\ssc^+$-sections.
A section $s\in \Gamma(p)$ satisfying $s(x)\in W(1)$ for all $x$ so that $s:X\rightarrow W(1)$ is sc-smooth is called a {\bf $\ssc^+$-section}.
The space of $\ssc^+$-sections is denoted by $\Gamma^+(p)$. In some sense $\ssc^+$-sections are compact perturbations, since the inclusion map
$W(1)\rightarrow W(0)$ is fiber-wise compact.

Finally we need the notion of a regularizing section.
\begin{defn}
Let $p:W\rightarrow X$ be a strong bundle over the M-polyfold $X$ and $f$ an sc-smooth section. We say that $f$ is {\bf regularizing} provided for a point $x\in X$ the assertion $f(x)\in W_{m,m+1}$ implies that $x\in X_{m+1}$.
\end{defn}

Note that for a regularizing section $f$ a solution $x$ of $f(x)=0$ belongs necessarily to $X_\infty$. If $f$ is regularizing and $s\in\Gamma^+(p)$, then $f+s$ is regularizing.

\section{Implicit Function Theorems}
The basic fact about the usual Fr\'echet differentiability is that if $f:U\rightarrow F$ is a smooth (in the usual sense) map between an open neighborhood
$U$ of $0$ in a Banach space $E$ and a Banach space $F$  and satisfies $f(0)=0$, then we can describe the solution set of $f=0$ near $0$ by an implicit function theorem
provided $df(0)$ is surjective and the kernel of $df(0)$ splits, i.e has a topological linear complement. So smoothness and some properties of the linearized operator at a solution give us always
qualitative knowledge about the solution set near $0$.  If on the other hand $f:U\rightarrow F$ is only sc-smooth and $Df(0)$ is surjective
and its kernel has an sc-complement we cannot conclude much about the solution space near $0$. However, as we shall see there is a large class
of sc-smooth maps for which a form of the implicit function theorem holds. In applications the class is large enough to explain gluing constructions (\`a la Taubes and Floer)
as   smooth implicit function theorems in the sc-world.

Our main aim is to define a suitable notion of Fredholm section of a polyfold bundle. One of the issues which has to be addressed at some point is the fact that the
spaces we are dealing with have locally varying dimensions. Though it might sound as a major issue it will turn out that there is a simple way to deal with  these type of problems. In fact, it is a crucial observation, that in applications base and fiber dimension change coherently. The sc-formalism incorporates this with a minimum amount of technicalities. One should remark that our presentation is slightly more general than that given in \cite{HWZ2}, and  stream-lines the presentation.

\subsection{A Special Class of Sc-Smooth Germs}
Let us begin with some notation. As usual $E$ is an sc-Banach space and $C\subset E$ a partial quadrant. We denote by $C_i$  the intersection $E_i\cap C$.  We shall write $\mathcal{O}(C,0)$ for an unspecified nested sequence
$U_0\supset U_1\supset U_2\supset \cdots$, where all the $U_i$ are relatively open neighborhoods of $0\in C_i\subset E_i$. Note that this differs from previous notation where $U_i=E_i\cap U$. When we are dealing with germs we always have the new definition in mind. A {\bf sc-smooth germ}
$$
f:{\mathcal O}(C,0)\rightarrow F
$$
is a map defined on $U_0$  so that for points $x\in U_1$ the tangent map
$Tf:U_1\oplus E_0\rightarrow TF$ is defined which again is a germ
$$
Tf:{\mathcal O}(TC,0)\rightarrow TF.
$$

We introduce a {\bf basic class} $\mathfrak{C}_{basic}$ of germs of maps as follows.
\begin{defn}
An element in $\mathfrak{C}_{basic}$ is an sc-smooth germ
$$
f:\mathcal{O}(([0,\infty)^k\times{\mathbb R}^{n-k})\oplus W,0)\rightarrow ({\mathbb R}^N\oplus W,0)
$$
for suitable $n$, $N$ and $0\leq k\leq n$, so that the following holds. If $P:{\mathbb R}^N\oplus W\rightarrow W$ is the projection, then $P\circ f$ has the form
$$
P\circ f(r,w)=w-B(r,w)
$$
for $(r,w)\in U_0\subset ([0,\infty)^k\times {\mathbb R}^{n-k})\oplus W$. Moreover, for every $\varepsilon>0$ and $m\in {\mathbb N}$
we have
$$
\norm{B(r,w)-B(r,w')}_m\leq \varepsilon\cdot \norm{w-w'}_m
$$
for all $(r,w), (r,w')\in U_m$ close enough to $(0,0)$ on level $m$.
\end{defn}
In \cite{HWZ2} the class of basic germs was slightly more general
in the sense that it was not required that $f(0)=0$ in its definition. However, all important results were then proved under the additional assumption that $f(0)=0$.
\begin{rem}
The thinking behind the above definition is implied by a look at the classical situation. Assume that $f:U\rightarrow F$ is a smooth (in the usual sense)
map defined on an open neighborhood $U$ of $0\in E$, where $F$ is a second Banach space. Suppose that $f(0)=0$ and $df(0)$ is a Fredholm map.
Denote by $K$ the kernel of $df(0)$ and by $C$ a complement of $df(0)E$.
If  $\dim(C)=N$ and
if $X$ is a complement of $K$, using that $df(0):X\rightarrow df(0)E$
is a topological linear isomorphism, we can define a topological linear isomorphism
$$
\phi:F=C\oplus df(0)E\rightarrow {\mathbb R}^N\oplus X
$$
in an obvious way, where $\phi$ has the property that
 $$
 \phi(0,h)=(0,(df(0)|X)^{-1}(h)).
 $$
 With $\dim(K)=n$ we can define
a topological linear isomorphism
$$
\varphi:E\rightarrow {\mathbb R}^n\oplus X
$$
which is the inverse of an isomorphism of the form
$$
\varphi:(a,x)\rightarrow \sigma(a)+x.
$$
Then
$$
E\oplus F\rightarrow ({\mathbb  R}^n\oplus X)\oplus ({\mathbb R}^N\oplus X)
$$
defined by
$$
(e,h)\rightarrow (\varphi(e),\phi(h))
$$
defines a Banach space bundle map which pushes $f$ forward to a map $g$
with $g(0,0)=(0,0)$ and $dg(0,0)(a,b)=(0,b)$. If $P:{\mathbb R}^N\oplus X\rightarrow X$ is the projection, then the map
$$
(a,x)\rightarrow Pg(a,x)=:x-B(a,x)
$$
has the property that for every $\varepsilon>0$  we have that
$$
\norm{B(a,x)-B(a,y)}\leq \varepsilon \cdot\norm{x-y}
$$
for all $(a,x)$ and $(a,y)$ close enough to $0$. So the fact that $f$ is smooth and has a Fredholm derivative gives a certain normal form after a suitable change of variables, which in our case  was quite trivial, namely
linear in the base and in the fiber.

In the applications to SFT and the other mentioned theories one can
bring the occurring nonlinear elliptic differential operators
even at bubbling-off points (modulo filling\footnote{This is a crucial concept in the polyfold theory and will be explained shortly.}) via sc-smooth coordinate changes in a similar form, in fact by quite conceptual methods, which are explained in
\cite{H2}  for Gromov-Witten theory, and in \cite{HWZ9} for the operators in SFT. It is important to note that if $f$ is sc-smooth so that $f(0)=0$ and $Df(0)$ is sc-Fredholm, it is generally not true that after a change of coordinates $f$ can be pushed forward to an element which belongs to $\mathfrak{C}_{basic}$.

As shown in \cite{HWZ2}, basic germs admit something like an infinitesimal
smooth implicit function theorem near $0$ (this is something intrinsic to sc-structures) which for certain maps can be `bound together' to a local implicit function theorem. To explain this, assume that $U\subset E$ is an open neighborhood of $0$
and $f:U\rightarrow F$ is an sc-smooth regularizing map satisfying $f(0)=0$ and $Df(0)$ is a surjective sc-Fredholm operator. Viewing $f$ as a section of $U\triangleleft F\rightarrow U$ assume that near every smooth point $x$ and
for a suitable $\ssc^+$-section with $s(x)=f(x)$ the germ $[f-s,x]$
is conjugated to a basic germ. Under these conditions there is a local implicit function theorem near $0$ which guarantees a local solution set
of dimension being the Fredholm index of $Df(0)$ at $0$ and in addition guarantees a natural manifold structure on this solution set.

The infinitesimal implicit function theorem refers to the following phenomena for basic germs. If $f\in \mathfrak{C}_{basic}$, then $Pf(a,w)=w-B(a,w)$,
where $B$ is a family of contractions on every level $m$ near $(0,0)$.
Hence, using Banach's fixed point theorem we find a germ $\delta_m$
solving $\delta_m(a)=B(a,\delta_m(a))$ on level $m$ for $a$ near $0$.
By uniqueness a solution on level $m$ also solves the problem on lower levels. This implies that we have a solution germ $a\rightarrow (a,\delta(a))$ of $Pf(a,w)=0$. The infinitesimal sc-smooth implicit function theorem gives the nontrivial fact that the germ
$$
\delta:{\mathcal O}([0,\infty)^k\times{\mathbb R}^{n-k},0)\rightarrow (W,0)
$$
is an sc-smooth germ.

In summary, as we shall discuss in more detail later, if we have a regularizing sc-smooth section which
around every smooth point is conjugated mod a suitable $\ssc^+$-section to a basic germ, then the 'infinitesimal' implicit function theorems around points $y$ near $x$, combine together to give a `local' implicit function theorem near a point $x$ where the linearization is surjective.
\end{rem}

\subsection{Sc-Fredholm Sections}
Assume next that $p:K\rightarrow O$ is a strong local bundle, i.e.
$(K,C\triangleleft F,E\triangleleft F)$  is a local strong bundle model.
Suppose $f$ is a germ of an  sc-smooth section near a smooth point $x\in O$, which we write as $f:{\mathcal O}(O,x)\rightarrow K$ and if there is no possibility of confusion we shall write $[f,x]$.
\begin{defn}
A {\bf filling} for the germ $[f,x]$ consists of the following data.
\begin{itemize}
\item[(1)] An  sc-smooth germ $\bar{f}:{\mathcal O}(C,x)\rightarrow F$.
\item[(2)] A choice of strong bundle retraction $R:U\triangleleft F\rightarrow U\triangleleft F$ such that $K$ is the image of $R$.
\end{itemize}
Viewing $f$ as a map $O\rightarrow F$
such that $\phi(y)f(y)=f(y)$, where $R(y,h)=(r(y),\phi(y)h)$ we assume that the data satisfies the following properties:
\begin{itemize}
\item[(1)] $\bar{f}(y)=f(y)$ for all $y\in O$ near $x$.
\item[(2)] $\bar{f}(y)=\phi(r(y))\bar{f}(y)$ for $y$ near $x$ in $U$ implies that $y\in O$.
\item[(3)] The linearisation of the map
$$
y\rightarrow (Id-\phi(r(y)))\bar{f}(y)
$$
at $x$ restricted to the $\ker(Dr(x))$ defines a linear topological isomorphism $\ker(Dr(x))\rightarrow \ker(\phi(x))$.
\end{itemize}
The germ $[f,x]$ is said to be {\bf fillable} provided there exists
a germ of strong bundle map $\Phi$, covering a local sc-diffeomorphism $\varphi$, so that the push-forward germ $[\Phi_\ast(f),\varphi(x)]$ has a filling. A {\bf filled version} of $[f,x]$ is an sc-smooth germ $[\bar{g},\bar{x}]$ obtained as a filling of a suitable push-forward.
\end{defn}

If $[f,x]$ has a filling $[\bar{f},x]$ the local study of $f(y)=0$ with $y\in O$ near $x$ is equivalent to the local study of $\bar{f}(y)=0$ where $y\in U$ close to $x$. Let us note that if $f(x)=0$ the linearisation $f'(x):T_xO\rightarrow K_x$ has the same kernel as $\bar{f}'(x):T_xU\rightarrow F_x$ and the cokernels are naturally isomorphic.

If $f$ is an sc-smooth section of a strong M-polyfold bundle $p:W\rightarrow X$ and $x$ is a smooth point we say that the germ $[f,x]$ admits a filled version, provided a local coordinate representation of $[f,x]$ admits a filled version as defined in the previous definition. We always may assume that the filled version has the form $g:{\mathcal O}(C,0)\rightarrow F$.

Now we come to the crucial definition.
\begin{defn}
We call the sc-smooth section $f$ of an  M-polyfold bundle $p:W\rightarrow X$ an  {\bf sc-Fredholm section}, provided $f$ is regularizing and around every smooth point $x$ the germ $[f,x]$ has a filled version $[g,0]$ so that for a suitable germ of $\ssc^+$-section $s$ with $s(0)=g(0)$ the germ $[g-s,0]$ is conjugated to an element in $\mathfrak{C}_{basic}$. We denote the collection of all sc-Fredholm sections of $p$ by ${\mathcal F}(p)$
\end{defn}

\begin{rem}
An sc-Fredholm section according to the above definition is slightly more general than the sc-Fredholm sections defined in \cite{HWZ2}. An additional advantage of the current definition is the stability result that
given an sc-Fredholm section for $p:W\rightarrow X$ and an $\ssc^+$-section $s$, then $f+s$ is an sc-Fredholm section for $p:W\rightarrow X$. With the version given in \cite{HWZ2} one can only conclude that $f+s$ is an sc-Fredholm section of $p^1:W^1\rightarrow X^1$. In applications the difference is only 'academic'. However, as far as a presentation is concerned this new version is more pleasant. This more general version will be implemented in \cite{H2,H3} and (conjecturally) even further improvements will be implemented in \cite{HWZ-polyfolds1,HWZ-polyfolds2}.
\end{rem}

The following stability result is crucial for the perturbation theory and rather tautological\footnote{In the set up of \cite{HWZ2} it was a nontrivial theorem. However,
some of the burden is now moved to the implicit function theorem.}.
\begin{thm}[Stability]
Let $p:W\rightarrow X$ a strong bundle over the M-polyfold $X$.
Then given $f\in {\mathcal F}(p)$ and $s\in \Gamma^+(p)$ we have
$f+s\in {\mathcal F}(p)$.
\end{thm}

Fredholm sections allow for an implicit function theorem.
We begin with the case without boundary.

\begin{thm}
Assume that $p:W\rightarrow X$ is a strong M-polyfold bundle, $f$ a sc-Fredholm section and $x$ a smooth point such that $f(x)=0$ and $f'(x):T_xX\rightarrow W_x$ is surjective. Then the solution set near $x$ carries in a natural
way the structure of a smooth manifold with dimension being the Fredholm index of $f'(x)$. In addition there exists an open neighborhood $V$ of $x$, so that for every $y\in V$ with $f(y)=0$
the linearisation $f'(y)$ is surjective. Moreover,  its kernel can be identified with the tangent spaces of the solution set at $y$.
\end{thm}

Next we start with the preparations for the boundary case.
Simple examples already show that the situation can be subtle.
Recall the Definition \ref{good-position} of a linear subspace being in good position
to a partial quadrant.

Assume that $X$ is an  M-polyfold (with boundary with corners) and
$x\in X$ a smooth point. The geometry of $X$ near $x$ is reflected in parts
by the degeneracy index $d:X\rightarrow {\mathbb N}$. If $x\in \partial X$,
then one might expect that the fact that $x$ is a boundary point
has a linearized version in $T_xX$ in the sense that there is a partial quadrant $C_x\subset T_xX$ so that the properties of $0\in T_xX$ in $C_x$ reflect those of $x\in X$. This is indeed true, see \cite{HWZ2}.
Namely consider all sc-smooth paths $\tau:[0,1]\rightarrow X$
starting at $x$ and take the closure of the set of all $\tau'(0)$, which will
define $C_x$. Now considering $C_x\subset T_xX$ we can talk about closed linear subspaces of $T_xX$ in good position to $C_x$. The degeneracy of $0\in C_x$ is precisely $d_X(x)$.

In the boundary case we have the following version of the implicit function theorem.
\begin{thm}
Assume that $p:W\rightarrow X$ is a strong M-polyfold bundle and $f$  an sc-Fredholm section. Assume that $x\in X$ is a smooth point satisfying $f(x)=0$ and $f'(x)$ is surjective. If the kernel $N$ of $f'(x)$ is in good position to $C_x$, the solution set $M$ of $f=0$ near $x$ carries in a natural way the structure of a smooth manifold with boundary with corners. For solutions $y\in M$ near $x$ the linearisations $f'(y)$ are also surjective and their kernels are also in good position to $C_y$.
\end{thm}

We refer the reader to \cite{HWZ2} for additional results and their proofs.

\subsection{Perturbations and Transversality}

Of course, like in the classical situation, a Fredholm section does not
need to be generic enough so that the previous results apply.
However, they can be brought into general position by perturbations
through $\ssc^+$-sections.

It is now important to understand what happens if we perturb a proper Fredholm section. Under which condition can we make sure that the perturbed
section is again proper?  In order to formulate some results we need
some auxiliary structures\footnote{In \cite{HWZ2} we introduced the notion of mixed convergence and that of an auxiliary norm. The definition of an auxiliary norm involves the notion of mixed convergence and the latter requires that the fibers of $W_{0,1}$ are reflexive. Though this is convenient in concrete applications an inspection of the proof shows that
for compactness assertions the notion of mixed convergence can be avoided
and the  definition of auxiliary norm can be simplified. We give the updated version in the present paper and comment on the necessary modifications. }.

An important metric concept is that of an auxiliary norm.
\begin{defn}
An {\bf auxiliary norm} $N$ for the strong M-polyfold bundle $p:W\rightarrow X$ consists of a continuous map $N:W_{0,1}\rightarrow [0,\infty)$ with the property  that the  restriction of $N$ to every fiber $(W_{0,1})_x$ is a complete norm.
\end{defn}
Using a continuous partition of unity on $X$ one can construct an auxiliary norm, see \cite{HWZ2} for the construction of auxiliary norms with additional useful properties.

We note the following useful lemma.
\begin{lem}
Assume that $p:W\rightarrow X$ is a strong M-polyfold bundle and
$K\subset X$ a compact subset. Given two auxiliary norms $N_1$ and $N_2$
there exists an open neighborhood $U$ of $K$ and a constant $c>0$ such that
$$
c\cdot N_1(h)\leq N_2(h)\leq c^{-1}\cdot N_1(h)
$$
for all $h\in W_{0,1}$ with $p(h)\in \overline{U}$.
\end{lem}

A modification of the  proof  from  \cite{HWZ2} gives the following theorem.
\begin{thm}\label{thk}
Assume that $p:W\rightarrow X$ is a strong polyfold bundle, $f$ an sc-Fredholm section of $p$, and $N$ an auxiliary norm.
Then for every smooth point $q\in X$  there exists an open neighborhood $U(q)\subset X$
so that every sequence $(x_k)\subset \overline{U(q)}$ satisfying
$\liminf_{k\rightarrow\infty} N(f(x_k))\leq 1$ has a convergent subsequence.
\end{thm}
We note the following consequence of the theorem.
\begin{rem}
The open neighborhood $U(q)$ guaranteed by Theorem \ref{thk} has the property that the set of points $\{x\in \overline{U(q)}\ |\ N(f(x))\leq 1\}$ is compact.
\end{rem}
Theorem \ref{thk} is an important input for the perturbation theory of proper Fredholm sections.
\begin{defn} An sc-smooth section $f$ is said to be {\bf proper} provided the solution set ${\mathcal S}:=\{x\in X\ |\ f(x)=0\}$ is compact in $X$.
\end{defn}
Using the regularizing property of a Fredholm section, it follows that
${\mathcal S}\subset X_\infty$. So one might ask whether ${\mathcal S}$ is a compact subset of $X_\infty$ with its finer topology. This is indeed true, but requires a proof which in fact uses strongly that a filled version has around every smooth point a germ conjugated to a basic germ, see \cite{HWZ2}.

\begin{prop}
If $f$ is a proper sc-Fredholm section of $p:W\rightarrow X$, then
the subset $\{x\in X\ |\ f(x)=0\}$ of $X_\infty$ is compact in $X_\infty$.
\end{prop}
Finally we arrive at an important result, which again is a consequence
of Theorem \ref{thk}, and which guarantees compactness for certain
perturbations of an sc-Fredholm section.
\begin{thm}\label{rot}
Let $p:W\rightarrow X$ be a strong bundle over the M-polyfold $X$, equipped with an auxiliary norm $N$. Finally let $f$ be a proper sc-Fredholm section. Then there exists an open neighborhood $U$ of ${\mathcal S}=\{x\in X\ |\ f(x)=0\}$ in $X$
so that for every $s\in \Gamma^+(p)$ with support in $U$, satisfying $N(s(x))\leq 1$ for all $x\in X$, the section $f+s$ is proper.
\end{thm}
Using the compactness of ${\mathcal S}$ and the local compactness the result
is quite obvious. We can cover ${\mathcal S}$ with finitely many open sets $U(q_1),\ldots ,U(q_n)$, where each set has the property guaranteed by the local compactness theorem. Let $U$ be their union and $s$ an $\ssc^+$-section with support in $U$ and $N(s(x))\leq 1$.
Next consider a sequence  $(x_k)$ of solutions of $f(x_k)+s(x_k)=0$. First of all $(x_k)$ has to be a subset of $U$. Further we see that $N(f(x_k))=N(s(x_k))\leq 1$ implying that $N(f(x_k))\leq 1$. After perhaps taking a subsequence we may assume that $x_k\in U(q_i)$ for some $i$.
The local compactness theorem now guarantees, that after taking a suitable subsequence, we may assume without loss of generality that $x_k\rightarrow x$.
From the regularizing property and $f(x_k)=-s(x_k)$ we conclude that $(x_k)\subset X_\infty$. Since $s(x_k)$ converges in $W_{0,1}$ to $s(x)$ we deduce from the normal form of $f$ that $x_k$ converges in $X_1$ which in turn gives that $s(x_k)$ converges in $W_{1,2}$. Boot-strapping our way up we see that $(x_k)$ converges in $X_\infty$. This,  in particular, implies that $f+s$ is proper sc-Fredholm section. Hence, if $s$ has support in $U$ and $N(s(x))\leq 1$, then $f+s$ is a proper sc-Fredholm section in ${\mathcal F}(p)$.

\begin{defn}
A pair $(U,N)$, where $U$ is an open neighborhood of the compact solution set ${\mathcal S}=\{x\in X\ |\ f(x)=0\}$ of a proper sc-Fredholms section,
and $N$ an auxiliary norm, which has the properties stated in Theorem \ref{rot}, is said to {\bf control compactness} for $f$.
\end{defn}

Having dealt with the compactness issue we turn our attention to transversality questions. If we want to associate to a proper sc-Fredholm section an invariant we need to be able to bring the section into
a sufficient general position\footnote{The wording has been chosen carefully here.
For example for Gromov-Witten theory one might read it as "general position".
In SFT there are many internal relationships between infinitely many different sc-Fredholm sections on polyfolds with boundaries to be preserved
and one cannot keep those, if one wants to bring the sections into a "general position". However, one can bring them into a position which still is
sufficiently general, keeping the internal relationships.}.

Here we only discuss the situation without boundary and refer the reader for the general discussion to \cite{HWZ2} and \cite{HWZ3}.
\begin{thm}
Let $p:W\rightarrow X$ be a strong bundle over the M-polyfold $X$, which is assumed not to have a boundary. Further, $f$ is a proper sc-Fredholm section and $N$ an auxiliary norm. Assume that $U$ is an open neighborhood  of $\{f=0\}$ so that $(U,N)$ controls compactness.
Then the following holds.
\begin{itemize}
\item[(1)] Given $\varepsilon \in (0,1)$ there exists a $\ssc^+$-section $s$ with support in $U$ and satisfying $N(s(x))<\varepsilon$ for all $x$ so that $f+s\in {\mathcal F}(p)$ has a surjective linearisation at every solution $x$ of $f(x)+s(x)=0$. In particular,  the solution set has in a natural way the structure of a compact smooth manifold.
\item[(2)] Given two such perturbations $s_0$ and $s_1$ as descried in (1) and viewing $f$ as a section of $W\rightarrow [0,1]\times X$, there exists an $\ssc^+$-section $s(t,x)$ of the latter, with support in $[0,1]\times U$ interpolating between $s_0$ and $s_1$, so that
    the sc-Fredhom section of $W\rightarrow [0,1]\times X$, defined
    by $(t,x)\rightarrow f(x)+s(t,x)$ is in transversal position to
    the zero section, and its solution set is a smooth compact manifold with boundary defining a cobordism between ${\mathcal S}_0$ and ${\mathcal S}_1$, where ${\mathcal S}_i=\{x\in X\ |\ f(x)+s_i(x)=0\}$.
\end{itemize}
\end{thm}
\begin{rem}
A slightly more general result can be proved, assuming that we start with
two different pairs $(U_i,N_i)$ controlling compactness and associated
$\ssc^+$-sections, which are as described in the theorem. Then one can construct first a connecting section for which compactness is assured, and then take a small perturbation to make it generic. The details
are basically similar as in the proof of the above theorem.
\end{rem}
\begin{rem}
One of the results concerning the boundary case, given in \cite{HWZ2}, states that
one can bring the solution set into general position to the boundary.
This,  in particular,  implies that the solution set is a manifold with boundary with corners, see \cite{HWZ2} for the M-polyfold version and \cite{HWZ3} for the polyfold version.
\end{rem}

We finish this section by introducing the notion of a differential form for M-polyfolds.

\begin{defn}
An {\bf sc-differential k-form } on the M-polyfold $X$ is an sc-smooth map
$\omega:\oplus_k TX\rightarrow {\mathbb R}$ which is separately linear in each argument and skew-symmetric. We denote the space of all forms (any degree) by $\Omega^\ast(X)$.
\end{defn}
One can define an exterior derivative $d$ and observes that $d\omega$ is a differential $(k+1)$-form on $X^1$. For practical purposes we need a space of forms, which is mapped via $d$ into itself. For this we observe that
the inclusion map $X^{i+1}\rightarrow X^i$ defines via pull-back a map
$\Omega^\ast(X^i)\rightarrow \Omega^\ast(X^{i+1})$. We denote by $\Omega^\ast_\infty(X)$ the direct limit and call its elements {\bf differential forms on $X_\infty$}. Then $d$ maps this space into itself and satisfies $d^2=0$. We can define therefore  deRham cohomology groups. It is trivial to develop an integration theory over submanifolds of M-polyfolds so that  Stokes' theorem holds. What is less trivial, see \cite{HWZ5},  is that this integration theory can be extended to much more complicated subsets which occur as solution sets
of the Fredholm theory in polyfolds, which we describe in the next section.

Let us recall that in the classical Fredholm theory we can consider orientable and oriented Fredholm sections. This is defined by means of an orientation of an associated line bundle, called the determinant bundle, see \cite{DK}. For these constructions it is important that the linearized operators depend continuously on the points where they were linearized. This is not true in our case.
However, as shown in \cite{HWZ8}, the construction can still be carried through. Assume that $E\rightarrow X$ is a strong  bundle over the M-polyfold $X$ and $f$ a Fredholm section. Give a smooth point $x\in X$ one can consider an $\ssc^+$-section $s$ with $f(x)=s(x)$ defined near $x$. Then the natural linearization of the section $f-s$ at $x$ is an sc-Fredholm section. Of course,
$(f-s)'(x)$ depends on $s$. Nevertheless, two different choices of $s$ result in two Fredholm operators differing by a linear $\ssc^+$-operator. In particular,  the space of all such linearizations is an affine space of sc-Fredholm operators.
Doing the above for all $x$ we obtain a bundle whose fibers
over $x$ are convex spaces of sc-Fredholm operators. If the determinant of one operator over a fixed $x$ has been oriented, all the others carry a natural orientation. An orientation for $f$ consists of a continuously varying (in $x$) orientation for all
these Fredholm operators\footnote{As pointed out there are some subtleties.}. An sc-Fredholm section $f$ of $E\rightarrow X$ equipped with an orientation is called an oriented Fredholm section. Observe that in the case of transversality the solution set $f^{-1}(0)$ of such a section has a natural orientation.

For example we have the following result, see \cite{HWZ2} for more results.
\begin{thm}
Let $E\rightarrow X$ be a strong bundle over the M-polyfold $X$
and $f$ a proper oriented Fredholm section. We assume that $X$ admits sc-smooth partition of unity\footnote{The sc-smooth partition of unity is needed to guarantee a sufficiently supply of $\ssc^+$-sections. The is for example guaranteed if $E_0$ is a Hilbert space.} and $\partial X=\emptyset$.
Then there exists a linear map
$$
\Phi_f:H_{dR}(X)\rightarrow {\mathbb R}
$$
characterized by the fact that for a sufficiently small generic perturbation $s$ (defined via a pair $(U,N)$) we have
$$
\Phi_f([\omega])=\int_{(f+s)^{-1}(0)} \omega.
$$
\end{thm}
The integral is defined to be $0$ if the dimensions do not match.

\section{Fredholm Theory in Polyfolds}
In this section we generalize the previous discussion to a level which is needed for the applications we have in mind.  We illustrate the ideas  by highlighting some of the issues in the  construction of Gromov-Witten invariants.
Full details are given in \cite{H2} and \cite{HWZ6}.

\subsection{Motivation via Gromov-Witten}

The following discussion is a slight modification from \cite{H2}.
We shall describe first a second countable paracompact space
$Z$, which is the space of stable curves (not necessarily pseudoholomorphic)
in a symplectic manifold, where the underlying maps have a certain
Sobolev regularity. As it turns out  this space carries `some kind of' sc-smooth structure. However, the situation is slightly more complicated, since the elements in $Z$ are in fact isomorphism classes of objects. In order to introduce  a smooth structure, one has to keep track of self-isomorphisms of objects which makes it necessary to generalize the discussion from the previous section.

Let $(Q,\omega)$ be a compact  symplectic manifold  without boundary. We  consider maps defined on
Riemann surfaces with images in  $(Q,\omega)$ having various
regularity properties. We shall write
$$
u:{\mathcal O}(S,x)\rightarrow Q
$$
for a mapping germ defined on a Riemann
surface $S$ near $x$.
\begin{defn}
Let $m\geq 2$ be an integer and $\varepsilon>0$.  We say a germ of
continuous map $u:{\mathcal O}(S,x)\rightarrow Q$ is of class
$(m,\varepsilon)$ at the point $x$ if for a smooth chart
$\phi:U(u(0))\rightarrow {\mathbb R}^{2n}$ mapping $u(0)$ to $0$ and
holomorphic polar coordinates $ \sigma:[0,\infty)\times
S^1\rightarrow S\setminus\{x\}$ around $x$, the map
$$
v(s,t)=\phi\circ u\circ \sigma(s,t)
$$
which is defined for $s$ large, has partial derivatives up to order
$m$, which weighted by $e^{\varepsilon s}$ belong to
$L^2([s_0,\infty)\times S^1,{\mathbb R}^{2n})$ if $s_0$ is sufficiently
large. We say the germ is of class $m$ around a point $z\in S$ provided
$u$ is of class $H^m_{loc}$ near $z$.
\end{defn}
 We observe that the above
definition does not depend on the choices involved, like charts and
holomorphic polar coordinates.

We consider now tuples $\alpha=(S,j,M,D,u)$, where $(S,j,M,D)$ is a
noded Riemann surface with ordered marked points $M$ and nodal pairs $D$,
and $u:S\rightarrow W$ is a continuous map.
\begin{defn}
A {\bf noded Riemann surface\index{Noded Riemann surface} with marked points}\index{Marked points} is a tuple $(S,j,M,D)$, where
$(S,j)$ is a closed Riemann surface, $M\subset S$ a finite collection of {\bf ordered marked points}, and $D$ is a finite collection of un-ordered pairs $\{x,y\}$ of points in $S$, called {\bf nodal pairs}, so that $x\neq y$ and two pairs which intersect are identical. The union of all $\{x,y\}$, denoted by $|D|$ is disjoint from $M$. We call $D$ the set of nodal pairs and $|D|$ the set of nodal points.
\end{defn}
The Riemann surface $S$ might consist of different connected components $C$. We call $C$ a {\bf domain component} of $S$. The special points\index{Special points} on $C$ are the points in $C\cap (M\cup|D|)$. We say that $(S,j,M,D)$ is {\bf connected}, provided the topological space $\bar{S}$ obtained by identifying $x=y$ in the nodal pairs $\{x,y\}\in D$ is connected.
With our terminology it is  possible that  $(S,j,M,D)$ is  connected  but on the other hand  $S$ may have several connected components, i.e. its domain components.

Next we describe the tuples $\alpha$ in more detail.

\begin{defn} We say that $\alpha=(S,j,M,D,u)$ is of {\bf class} $({\bf m},\boldsymbol{\delta})$
provided the following holds, here $m\geq 2$ and $\delta>0$.
\begin{itemize}
\item[(0)] The underlying topological space obtained by identifying
the two points in any nodal pair is connected\footnote{For certain applications one would allow non-connected $(S,j,M,D)$.}.
\item[(1)] The map $u$ is  of class $(m,\delta)$ around the points in $|D|$ and
of class $m$ around all other points in $S$ \footnote{For certain
applications it is useful to require  the map $u$ around the points in $M$ to be of class $(m,\delta)$ as well.
This, however, requires only minor modifications.}.
\item[(2)] $u(x)=u(y)$ for every nodal pair $\{x,y\}\in D$.
\item[(3)] If a  domain component $C$ of $S$ has genus $g_C$ and
$n_C$ special points so that $2\cdot g_C +n_C\leq 2$, then $\int_C u^\ast\omega >0$.
\end{itemize}
We call two such tuples $(S,j,M,D,u)$ and $(S',j',M',D',u')$
{\bf equivalent}\index{Equivalence of stable maps} if there exists a biholomorphic map $
\phi:(S,j,M,D)\rightarrow (S',j',M',D') $ with $$u'\circ\phi=u.$$ An
equivalence class is called a (connected) {\bf stable curve of class}
$({\bf m},\boldsymbol{\delta})$\index{Stable curve of class $(m,\delta)$}. An $\alpha$ is called a {\bf stable map}. Hence a stable curve is an equivalence class of stable maps. Requirement  (3) is called the {\bf stability condition}.
\end{defn}

Next we introduce the space $Z$ which later will be equipped with
the polyfold structure.
\begin{defn}
Fix a $\delta_0\in (0,2\pi)$. The collection of all equivalence
classes $[\alpha]$ of tuples $\alpha$ of class $(3,\delta_0)$ is called the
{\bf space of stable curves} into $Q$ of class $(3,\delta_0)$ and is denoted
by $Z$, or by $Z^{3,\delta_0}(Q,\omega)$.
\end{defn}
The set $Z$ has a natural topology. More precisely:
\begin{thm}\label{th-top}
Given $\delta_0\in (0,2\pi)$ the space $Z^{3,\delta_0}(Q,\omega)$
has a natural second countable paracompact topology.
\end{thm}
We refer the reader to \cite{H2,HWZ6} for more detail.

Assume next that a compatible almost complex structure $J$ has been fixed for $Q$.
We consider tuples $(S,j,M,D,u,\xi)$, where $(S,j,M,D,u)$ is as just described and $\xi(z):T_zS\rightarrow T_{u(z)}Q$ is complex anti-linear
for the given structures $j$ and $J$. Further we assume that $\xi$ has Sobolev regularity $H^2$ away from the nodal points. At the nodal points
we assume it to be of class $(2,\delta_0)$. To make this precise,
pick a nodal point $x$ and take positive holomorphic polar coordinates around $x$, say $\sigma(s,t)$ with $x=\lim_{s\rightarrow \infty} \sigma(s,t)$.
Then take a smooth chart $\phi$ around $u(x)$ with $\varphi(u(x))=0$.
Finally consider
$$
(s,t)\rightarrow pr_2\circ T\varphi(u(\sigma(s,t)))\xi(\sigma(s,t))(\frac{\partial\sigma(s,t)}{\partial s}),
$$
which we assume for large $s_0$ to be in $H^{2,\delta_0}([s_0,\infty)\times S^1,{\mathbb R}^{2n})$. The definition of the decay property does not depend on the choice of $\sigma$ and $\varphi$. The previously defined notion of equivalence for stable maps extends to these tuples. Denote the collection of all these by $W=W^{2,\delta_0}$. One can equip this space with a natural second countable paracompact topology so that we have, in particular, the natural surjective and continuous map
$$
p:W\rightarrow Z:[S,j,M,D,u,\xi]\rightarrow [S,j,M,D,u].
  $$
  As we shall see $p$ can be equipped with the structure of a so-called strong polyfold bundle. Further, the natural section $\bar{\partial}_J$ defined by
$$
\bar{\partial}_J([S,j,M,D,u])=[S,j,M,D,u,\bar{\partial}_{J,j}(u)]
$$
defines a (polyfold-) Fredholm section, which on every connected component of $Z$ is proper. The details are provided in \cite{H2,HWZ6}. In the next subsection
we shall describe the polyfold theory. We start however with some motivation coming from the current example.\\

\noindent {\bf Motivation for the polyfold construction:} Consider the category $\mathfrak{C}$ whose objects are the tuples $\alpha=
(S,j,M,D,u)$. This is a very large category, which is not even a set. Next we define morphisms as follows. A morphism $\Phi:\alpha\rightarrow \alpha'$ is a tuple $(\alpha,\phi,\alpha')$, where $\phi:(S,j,M,D)\rightarrow (S',j',M',D')$
is a biholomorphic map preserving the ordered marked points $M$ and the nodal pairs, so that in addition $u'\circ \phi=u$. If we call two objects equivalent if there is a morphism between them we obtain the orbit space $|\mathfrak{C}|=\mathfrak{C}/\sim$. The latter is by our definition precisely $Z$, which is a set.

The good thing about $\mathfrak{C}$ compared to $Z$ is that we keep track of the symmetries, the bad thing (among others) is that $\mathfrak{C}$ is too big.
The first idea is now to `weed' out a lot of objects so that we obtain a small category, i.e. the objects and morphisms form a set. Of course, we insist that the associated orbit space is still $Z$. The `weeding out' has to be done cleverly. Namely we want to be able
to define a topology on the set of objects which gives us a notion about closeness of two objects or two morphisms (Note that we talk about closeness of two objects and not about the closeness of two equivalence classes!). If that has been achieved we even try
more, namely to define M-polyfold structures on the topological spaces of objects and morphisms so that the natural operations in the category become smooth.
This will lead to the notion of an ep-groupoid (or polyfold groupoid), which can be viewed as a generalization of an atlas of an  M-polyfold. The above, suitably modified
can, of course, be also used to motivate orbifolds and \'etale proper Lie groupoids. In fact, ep-groupoids can be viewed as the M-polyfold generalization of \'etale proper Lie groupoids.

\subsection{Polyfold Groupoids and Polyfolds}
The following material is discussed in detail in \cite{H2,HWZ3}. The reader should note that in some sense we generalize the notion of an \'etale proper Lie groupoid to the differential geometry based on local sc-models. As such we generalize ideas described in \cite{Mj}, which go back to ideas presented in
\cite{hae1,hae2,hae3}.

A {\bf groupoid}  ${\mathcal G}$ is a small category, where every morphism
is invertible. Recall that a category is small, provided the class of objects as well as morphisms is a set. In the future we shall denote the groupoid as well as its object set by $G$ and the morphism set by ${\bf G}$. From the category structure we have five structure maps. The {\bf source-} and  {\bf target maps} $t,s:{\bf G}\rightarrow G$ associating to a morphism its source and target, respectively. Then there is the {\bf inversion map} $i:{\bf G}\rightarrow {\bf G}$, the {\bf $\boldsymbol{1}$-map} $G\rightarrow {\bf G}:x\rightarrow 1_x$ and the {\bf multiplication map} $m:{\bf G}{{_s}\times_t}{\bf G}\rightarrow {\bf G};(\phi,\psi)\rightarrow \phi\circ\psi$.

In the following by an  M-polyfold  we mean  an  sc-smooth space which is build on neat local sc-models. Of course, a subcase is that where we do not have a boundary
with corners.

\begin{defn}
An  {\bf ep-groupoid}, also called {\bf polyfold groupoid},  consists of an groupoid $X$  where the object and morphism sets carry M-polyfold structures, so that the following holds.
\begin{itemize}
\item[(1)] The source and target maps are surjective local sc-diffeomorphisms. In particular,  the fibered product ${\bf X}{{_s}\times_t}{\bf X}$ has a natural M-polyfold structure\footnote{It can be shown that ${\bf X}{{_s}\times_t}{\bf X}$ has in a natural way an M-polyfold structure since $t$ and $s$ are local sc-diffeomorphisms.}.
\item[(2)] The structure maps $i:{\bf X}\rightarrow {\bf X}$, $u:X\rightarrow{\bf X}$ and $m:{\bf X}{{_s}\times_t}{\bf X}\rightarrow {\bf X}$ are sc-smooth.
\item[(3)] Every point $x\in X$ has an open neighborhood $U=U(x)$ so that $t:s^{-1}(\overline{U})\rightarrow X$ is proper\footnote{The usual definition of properness as defined in the context of Lie groupoids cannot be carried over to our possibly infinite-dimensional case.}. (Using that the inversion map is an  sc-diffeomorphism one can interchange the role of $s$ and $t$ in this definition.)
    \end{itemize}
    \end{defn}
    Given an ep-groupoid $X$ we can define its {\bf orbit space} $|X|$  by calling two points $x$ and $y$ in $X$ equivalent providing there exists a morphism $\phi:x\rightarrow y$. We equip $|X|$ with the quotient topology. Clearly $|X|$ has a filtration $|X|_m:=|X_m|$ by nested topological spaces.

    There is an important structural statement about the local geometry of
    an ep-groupoid which is proved in \cite{HWZ3} and is a generalization of a similar fact in the theory of Lie groupoids in \cite{Mj}.
 \begin{thm}
Given an ep-groupoid $X$ and $x\in X$, there exist  an open
neighborhood $U\subset X$ of $x$,  a group homomorphism
$$
\varphi:G_{x}\rightarrow \hbox{Diff}_{\ssc}(U), \quad g\mapsto
\varphi_g,
$$
and an  sc-smooth map
$$
\Gamma:G_x\times U\rightarrow {\bf X}
$$
having the following properties.
\begin{itemize}
\item[$\bullet$] $\Gamma(g,x)=g$.
\item[$\bullet$] $s(\Gamma(g,y))=y$ and $t(\Gamma(g,y))=\varphi_g(y)$.
\item[$\bullet$] If $h:y\rightarrow z$  is a morphism with $y,z\in U$,  then there exists a
unique $g\in G_x$ with $\Gamma(g,y)=h$.
\end{itemize}
\end{thm}
In particular, every morphism between points in $U$ belongs to the
image of the map $\Gamma$.  We call $\varphi:G_x\rightarrow
\hbox{Diff}_{\ssc}(U)$ a {\bf natural representation} of the
stabilizer $G_x$.

To continue the discussion we need the following definition.
\begin{defn}
A functor $F:X\rightarrow Y$\index{$\ssc^k$-functor} between two ep-groupoids is $\ssc^k$, \index{Sc-smooth functors} provided the induced map between the object and morphism M-polyfold is $\ssc^k$.
\end{defn}

Let us observe that  an $\ssc^k$-functors induces an  $\ssc^0$-map between the orbit spaces $|F|:|X|\rightarrow |Y|$.

Let us denote for $x\in G$ by ${\bf G}(x)$ the isotropy group of $x$. It consists of all morphisms $\phi:x\rightarrow x$. Observe that for an ep-groupoid $G$ every isotropy group is finite. This is an immediate consequence of the properness condition and the fact that the source and target maps are local sc-diffeomorphisms.
\begin{defn}
Let $X$ and $Y$ be two polyfold groupoids. An {\bf equivalence} $F:X\rightarrow Y$ is a functor satisfying the following
\begin{itemize}
\item[(1)] $F$ is a local sc-diffeomorphism on objects.
\item[(2)] $F$ induces a sc-homeomorhism $|F|:|X|\rightarrow |Y|$.
\item[(3)] For every $x\in G$ the functor $F$ induces a bijection between the isotropy groups ${\bf G}(x)$ and ${\bf G}(F(x))$.
    \end{itemize}
    \end{defn}

 \begin{rem}
 Equivalences are in general not invertible as functors. However, this picture will improve when we pass to some category with the same objects in which these equivalences can be inverted.
At the moment we consider the category which has as objects the polyfold groupoids  and the sc-smooth functors as morphisms. If we view an ep-groupoid as a generalization of an atlas we would like for example that a refinement of an atlas is an isomorphic object.
 An equivalence can be viewed as the generalization of such a refinement.
 However, as we have seen it is not invertible in general. Using a standard construction from category theory, \cite{GZ}, it is always possible to invert a given family of arrows, by keeping the objects, and by only making a minimal amount (depending on the original family of arrows) of "structural damage" to the category. The concrete realization of this procedure
 in our category has a simple description, which we will describe later and where we follow the ideas from the Lie groupoid case, see \cite{Mj}.
 See \cite{GZ} for the general theory.
 \end{rem}

 We want to relate two polyfold groupoids and introduce the notion of a common refinement, which in some sense generalizes the idea that two (compatible) atlases on a manifold have a common refinement.

\begin{defn}
Let $X$ and $X'$ be two polyfold groupoids. A {\bf  common refinement} consists of a third polyfold groupoid $X''$ and equivalences $F:X''\rightarrow X$ and $F':X''\rightarrow X'$.
\end{defn}

An important observation is that this defines an equivalence relation.
\begin{prop}
The notion of having a common refinement is an equivalence relation.
\end{prop}

We need the notion of natural equivalence in our category of polyfold groupoids.
\begin{defn}
Two sc-smooth functors $F,G:X\rightarrow Y$ between the same pair of polyfold groupoids are called {\bf naturally equivalent} provided there exists an  sc-smooth map
$\tau:X\rightarrow {\bf Y}$ associating to an object $x\in X$ a morphism $\tau(x):F(x)\rightarrow G(x)$ which is natural in the sense that for every $h:x\rightarrow x'$ in ${\bf X}$ the identity
$$
\tau(x')\circ F(h)=G(h)\circ \tau(x)
$$
holds. The map $\tau$ is called a {\bf natural transformation}.
\end{defn}
 Let us observe that two naturally equivalent functors induce the same maps $|F|$ and $|G|$ between the orbit spaces. We also observe that natural equivalence is an equivalence relation.

As we already mentioned before, we would like to construct a new category
in which equivalences are invertible. In order to do so we consider diagrams\index{Diagrams of sc-smooth functors}
$$
d: X\xleftarrow{F} A\xrightarrow{\Phi}Y,
$$
where $X$, $A$ and $Y$ are polyfold groupoids, $F$ is an equivalence and $\Phi$ an sc-smooth functor. Such a diagram will, as we shall see, be a representative
of a morphism between $X$ and $Y$. Let us call such a diagram $d$ a diagram from $X$ to $Y$. We observe that such a diagram induces an $\ssc^0$-map $|d|$ defined by
$$
|d|=|\Phi|\circ |F|^{-1}:|X|\rightarrow |Y|.
$$
In a first step we define the notion of a refinement of such a diagram.
\begin{defn}
Assume that $d:X\xleftarrow{F} A\xrightarrow{\Phi}Y$ is a diagram from $X$ to $Y$. A diagram $d':X\xleftarrow{F'} A'\xrightarrow{\Phi'}Y$  from $X$ to $Y$ will be called a {\bf refinement} of $d$ if there exists an equivalence $H:A'\rightarrow A$ such that $F\circ H$ and $F'$ are naturally equivalent and $\Phi\circ H$ and $\Phi'$ are naturally equivalent. Let us write $H:d'\rightarrow d$ if $d'$ refines $d$
via the equivalence $H':A'\rightarrow A$.
\end{defn}
Here are two remarks:
\begin{rem} \mbox{}\\
(1) Trivially, if $d'$ refines $d$, the induced maps $|d|$ and $|d'|$ are the same.\\
\noindent (2) If $d'$ refines $d$, and $d''$ refines $d'$, then $d''$ refines $d$.
\end{rem}

It is clear what it means that two diagrams $d$ and $d'$ from $X$ to $Y$ have a common refinement. The crucial observation is the following proposition.
\begin{prop}
Assume that $d$, $d'$ and $d''$ are diagrams from $X$ to $Y$. Assume that
$d$ and $d'$, and $d'$ and $d''$ have common refinements. Then $d$ and $d''$ have a common refinement.
\end{prop}
As a consequence of the proposition we see that saying the two diagrams from $X$ to $Y$ have a common refinement defines an equivalence relation.
\begin{defn}
The equivalence class $[d]$ of a diagram from $X$ to $Y$ is called a {\bf generalized map}.
\end{defn}
By a previous discussion, we know that there is a well-defined map
$|[d]|:=|d|:|X|\rightarrow |Y|$ associated to such an equivalence class.

Let us note that if $[d]:X\rightarrow Y$ and $[e]:Y\rightarrow Z$ are generalized maps, there is a well-defined composition $[e]\circ [d]:=[f]$, where $f$ is suitably defined. The construction is nontrivial, see \cite{HWZ3}, and needs the construction of weakly fibered products in the category of ep-groupoids. Of course, the construction is known in the Lie groupoid context.

At this point, the category of interest to  us, is the category where
the objects are polyfold groupoids and the morphisms are generalized maps, i.e.
equivalence classes of diagrams. The original category with the same objects, but the sc-smooth functors as morphisms, maps  into this new category
as the identity on objects, and on morphisms, by mapping an sc-smooth functor
$\Phi:X\rightarrow Y$ to the equivalence class of the diagram
$$
X\xleftarrow{Id} X\xrightarrow{\Phi}Y.
$$
If we start with an equivalence $F:X\rightarrow Y$ the associated diagram $d_F$
has an equivalence class $[d_F]$ which is inverted by the equivalence class of the diagram $Y\xleftarrow{F} X\xrightarrow{Id} X$. In particular, any diagram
$X\xleftarrow{F} Y\xrightarrow{F} X$, where $F$ is an equivalence represents the identity morphism of $X$.
\begin{defn}
A generalized map $[d]:X\rightarrow Y$ is said to be {\bf strongly invertible} (s-invertible for short) if there exists a representative $d$ of the form $X\xleftarrow{F} A\xrightarrow{G}Y$
with $F$ and $G$ being equivalences.
\end{defn}

Now we can define the notion of a polyfold structure.
\begin{defn}
A {\bf polyfold structure} on a  metrizable  space $Z$ is a pair $(X,\beta)$, where $X$ is a polyfold groupoid and $\beta:|X|\rightarrow Z$ a homeomorphism.  Two polyfold structures $(X,\beta)$ and $(X',\beta')$ on $Z$ are equivalent provided there exists an s-invertible generalized map $[d]:X\rightarrow X'$ with $\beta'\circ |d|=\beta$. A metrizable space equipped with an equivalence class of polyfold structures is called  a {\bf polyfold}.
\end{defn}
Given a polyfold $Z$ and $(X,\beta)$ being one of the defining structures, we call $X$ a {\bf model} for $Z$.

Let us observe that a polyfold $Z$ has a degeneracy map $d:Z\rightarrow {\mathbb N}$  defined by $d(z):=d_X(\beta^{-1}(z))$, where $X$ is a polyfold groupoid so that $(X,\beta)$ defines the polyfold structure. This is independent of the representative since equivalences preserve the degeneracy index being local sc-diffeomorphisms.

\begin{defn} Let $Z$ and $Z'$ be two polyfolds with polyfold structures defined by $(X,\beta)$ and $(X',\beta')$, respectively. An sc-smooth map  $f:Z\rightarrow Z'$,
is an equivalence class of pairs $[(f,[d])]$, where $[d]:X\rightarrow X'$ is a generalized map so that $\beta'\circ |d| = f\circ\beta$. It is clear what the notion of equivalence is.
\end{defn}

So considering  an sc-smooth map $f:Z\rightarrow Z'$ is equivalent to considering
 an  sc-smooth functor $\Phi:A\rightarrow X'$ for suitable ep-groupoids $A$ and $X'$. Indeed, if $(X,\beta)$ is the model for $Z$ and $(X',\beta')$ for $Z'$
and $d:X\xleftarrow{F}A\xrightarrow{\Phi} X'$ satisfies  $\beta'\circ |d|=f\circ \beta$, then $(A,\beta\circ |F|)$ is an equivalent model for $Z$ and with respect to this model $f$ is represented by $\Phi:A\rightarrow X'$.

There is a parallel discussion for strong bundles over ep-groupoids or polyfolds.
We refer the reader to \cite{HWZ3} for details. A strong bundle over an ep-groupoid is a strong bundle $E\rightarrow X$ over the object M-polyfold
together with a certain strong bundle map $\mu$. In order to explain this further, we observe that the pull-back of $p:E\rightarrow X$ via $s:{\bf X}\rightarrow X$ defines the strong bundle ${\bf X}{{_s}\times_p}E\rightarrow{\bf X}$. Then $\mu$ is a strong bundle map
$$
{\bf X}{{_s}\times_p}E\rightarrow E
$$
covering $t:{\bf X}\rightarrow X$, which is a linear isomorphism in the fibers, so that writing $\phi\cdot e:=\mu(\phi,e)$ we have
$$
1_x\cdot e =e
$$
for $e\in E_x$ and
$$
\psi\cdot (\phi\cdot e)=(\psi\circ\phi)\cdot e.
$$
In  other words the morphisms in $X$ are lifted to the fiber.
 One can  construct from the data $E\rightarrow X$ and $\mu$ a morphism set ${\bf E}$ for $E$, which turns $E$ into an ep-groupoid.
In fact ${\bf E}$ will have the structure of a strong bundle over ${\bf X}$.
Then $p:E\rightarrow X$ will be an sc-smooth functor, see \cite{HWZ3} for details.

Given two metrizable spaces and a surjective map $p:W\rightarrow Z$,
a structure of a strong polyfold bundle is given by a strong bundle over an ep-groupoid together with a homeomorphism of the orbits spaces with $W$ and $Z$ preserving $p$.
As in the polyfold case one defines similarly a notion of equivalence.
Then $p:W\rightarrow Z$ equipped with an equivalence class of strong bundle structures
is called  strong polyfold bundle.

In the case of our example, the Gromov-Witten theory, we have the following.
Recall the "bundle"  $p:W\rightarrow Z$.

\begin{thm}
Given a strictly increasing sequence of weights
$(\delta_k)$ starting with $0< \delta _0<2\pi$ and staying below $2\pi$, there exists  a strong polyfold bundle structure on $p:W\rightarrow Z$,  so that $\bar{\partial}_J$ is an sc-smooth section.
\end{thm}

For a proof see \cite{HWZ6}.

We need to generalize the notion of a submanifold of an  M-polyfold to ep-groupoids and polyfolds. There are, in fact, many possibilities, for example
one might consider subsets of a polyfold which have a natural orbifold structure. However, our guiding principle is Fredholm theory and the associated transversality questions, which seriously restricts the possibilities. Even in finite dimension it is not difficult to construct an orbifold bundle for which all smooth sections will not be transversal to the zero-section. The functoriality of a description in the model
gives obstructions to transversality. We can, however, achieve transversality by breaking symmetries. Of course, this is bad, since it destroys the functorial structure. The best fix for this is to consider locally a family
of problems, so that the family is invariant under the symmetries.
This is the road we take leading us to multi-valued perturbations.
It generalizes some of the work in a more standard setting of global group actions on Hilbert manifolds in \cite{CRS}, to the categorical setting of ep-groupoids and polyfolds. We shall introduce the notion of a branched ep-subgroupoid. It is in some sense the most general type of smooth subset
which we can get from a transversally perturbed Fredholm functor using multisections. In certain situations, i.e. having more knowledge,
we can sometimes expect generic solution sets to have more structure.

For the following, view the set of non-negative rational numbers as a category with the only morphisms being the identities. We recall a  definition  from \cite{HWZ5}.

\begin{defn}
A {\bf branched ep-subgroupoid} of dimension $n$ of the ep-groupoid $X$ is a functor $\Theta:X\rightarrow {\mathbb Q}^+$ with the following properties:
\begin{itemize}
\item[(1)] The support of $\Theta$, defined by $\hbox{supp}(\Theta)=\{x\in X\ |\ \Theta(x)>0\}$ is contained in $X_\infty$.
    \item[(2)] Every point $x$ in the support of $\Theta$ is contained in an open neighborhood $U=U(x)$ such that
        $$
        \hbox{supp}(\Theta)\cap U =\bigcup_{i\in I} M_i
        $$
where $I$ is a finite set and where the $M_i$ are finite-dimensional submanifolds of $X$ of dimension $n$.
\item[(3)] There exist positive rational numbers $\sigma_i$, $i\in I$, such that for $y\in U$ we have
    $$
    \Theta(y)=\sum_{\{i\in I\ |\ y\in M_i\}} \sigma_i.
    $$
    \item[(4)] The inclusion maps $M_i\rightarrow U$
    are proper.
    \item[(5)] There is a natural representation of the isotropy group $G_x$ on $U$.
        \end{itemize}
        \end{defn}
Given $\Theta$ one can introduce the notion of orientability and orientation. We refer the reader to \cite{HWZ5} for more details.

We can extend the above definition of a branched ep-subgroupoid to polyfolds, see \cite{HWZ5} for all details.
\begin{defn}
A {\bf branched suborbifold } $S$ of a polyfold $Z$ is a subset $S$ of $Z$  equipped with a weight function $w:S\rightarrow (0,\infty)\cap{\mathbb Q}^+$ together with an equivalence class of triples $(X,\beta,\Theta)$, where $X$ is an ep-groupoid, $(X,\beta)$  a model for $Z$, and $\Theta$ a branched ep-subgroupoid, so that the image of $\hbox{supp}(\Theta)$ under the map
$x\rightarrow \beta(|x|)$ is $S$ and  $\Theta(x)=w\circ \beta(|x|)$ for all $x\in \hbox{supp}(\Theta)$.
\end{defn}
The equivalence of tuples $(X,\beta,\Theta)$ is as usual defined via refinements and equivalences. Again one can define orientatibility and orientation. It is, of course, done via the overhead using the same notions
occurring in the ep-groupoid case.

In order to define invariants for polyfold Fredholm sections
it will be important to define the notion of a sc-differential form on an ep-groupoid. We already introduced this  notion for M-polyfolds. In the case of an ep-groupoid $X$, an sc-differential form is simply an sc-differential form, say $\omega$, on the object M-polyfold, so that for every morphism $\phi:x\rightarrow y$, where $x$ is at least on level $1$, we have
$$
\omega_y\circ \oplus_kT\phi =\omega_x.
$$
Again we can define $\Omega^\ast_\infty(X)$ as in the M-polyfold case, except that we also require the compatibility with morphisms. Then we have an exterior derivative and can define
a version of deRham cohomology
$$
H^\ast_{dR}(X)\ \hbox{and}\ \ H^\ast_{dR}(X,\partial X),
 $$
 see \cite{HWZ5} for details.

This can be extended to polyfold  $Z$ and we obtain the notion of an sc-differential form on $Z_\infty$ and associated deRham groups $$
H^\ast_{dR}(Z)\ \ \hbox{and}\ \  H^\ast_{dR}(Z,\partial Z).
$$

As shown in \cite{HWZ5} there is an appropriate  integration theory, which allows to integrate sc-differential forms
on an ep-groupoid over a branched ep-subgroupoid. The naturality of the integration makes it compatible with equivalences and therefore is defined for branched sub-orbifolds of polyfolds.

In the following $G_e$ denotes the effective part of the automorphism group $G$. To explain this call an automorphism $\phi:x\rightarrow x$ non-effective if the associated $t\circ s^{-1}$ is the identity. It is not difficult to show that the non-effective elements form a normal subgroup $G_{ne}$ so that
$G_e:=G/G_{ne}$ is a group called the effective part.

\begin{thm}[Canonical Measures]
Let $X$ be an ep-groupoid and assume that $\Theta:X\to \Q^+$ is an oriented branched
ep-subgroupoid of dimension $n$ whose orbit space $S=\abs{\supp \Theta}$  is compact and equipped with the weight function
$\vartheta:S\to \Q^+$ defined by
$$
\vartheta (|x|):=\Theta (x),\quad |x|\in S.
$$
Then there exists a map
$$
\Phi_{(S,\vartheta)}:\Omega^n_\infty(X)\rightarrow {\mathcal M}(S,{\mathcal
L}(S)),\quad \omega \mapsto    \mu_\omega^{(S,\theta)}
$$
which  associates to every  sc-differential $n$-form $\omega$ on $X_\infty$ a signed finite measure
$$\mu_{\omega}^{(S, \vartheta )}\equiv \mu_{\omega}$$
on  the canonical measure space $(S,{\mathcal L}(S))$\footnote{Recall that a smooth manifold admits in a natural way a $\sigma$-algebra of measurable sets, called the Lebesgue $\sigma$-algebra. This is a generalization.}. This map is uniquely characterized by the
following properties.
\begin{itemize}
\item[(1)] The map $\Phi_{(S,\theta)}$ is linear.
\item[(2)] If $\alpha=f\tau$ where $f\in\Omega^0_\infty(X)$ and $\tau\in\Omega^{n}_\infty(X)$, then
$$
\mu_\alpha(K)=\int_K
fd\mu_\tau
$$
for  every set $K\subset S$ in the $\sigma$-algebra ${\mathcal L}(S)$.
\item[(3)] Given a  point $x\in \supp \Theta$ and an oriented branching structure
${(M_i)}_{i\in I}$  with  the associated weights $(\sigma_i)_{i\in I}$ on the open neighborhood  $U$ of $x$, then  for every set  $K\in {\mathcal L}(S)$ contained in a compact subset of $\abs{\supp \Theta \cap U}$,  the $\mu_{\omega}$-measure of $K$ is given by the formula
$$
\mu_\omega (K)=\frac{1}{\sharp G_{e}}\sum_{i\in I}
\sigma_i\int_{K_i}\omega\vert M_i
$$
where $K_i\subset M_i$ is the preimage of $K$ under  the projection map
$M_i\to \abs{\supp \Theta \cap U}$ defined by $x\rightarrow |x|$.
\end{itemize}
\end{thm}

Here $\int_{K_i}\omega\vert M_i$ is the signed measure of the set
$K_i$ with respect to the Lebesgue signed measure associated with the smooth $n$-form $\omega\vert M_i$ on the finite dimensional manifold $M_i$.
A similar result holds for the boundary.

\begin{thm}[Canonical Boundary  Measures] Under the same assumptions
as above there exists a map
$$
\Phi_{(\partial S,\vartheta)}:\Omega^{n-1}_\infty(X)\rightarrow {\mathcal
M}(\partial S,{\mathcal L}(\partial S)), \quad \tau\mapsto
\mu_\tau^{(\partial S,\vartheta)}
$$
which assigns to every sc-differential $(n-1)$-form $\tau$ on $X_\infty$ a signed finite measure
$$\mu_{\tau}^{(\partial S, \vartheta)}\equiv \mu_{\tau}$$
on  the canonical measure  space
$(\partial  S, {\mathcal L} (\partial S))$.
 This map is uniquely characterized by the
following properties.
\begin{itemize}
\item[(1)] The map $\Phi_{(\partial S,\vartheta )} $ is linear.
\item[(2)] If $\alpha=f\tau$ where  $f\in\Omega^0_\infty(X)$ and $\tau\in\Omega^{n-1}_\infty(X)$, then every $K\in {\mathcal L}(\partial S)$
has the $\mu_{\alpha}$-measure
$$\mu_\alpha(K)=\int_K
fd\mu_\tau.$$
\item[(3)] Given a  point $x\in \supp \Theta \cap \partial X$  and an oriented branching structure  ${(M_i)}_{i\in I}$ with weights $(\sigma_i)_{i\in I}$ on the open neighborhood $U\subset X$ of $x$, then the measure of  $K\in {\mathcal L}(\partial S)$  contained in a compact subset of $\abs{\supp \Theta \cap U\cap \partial X}$ is given by the formula
$$
\mu_\tau(K)=\frac{1}{\sharp G_{e}}\sum_{i\in I}
\sigma_i\int_{K_i}\tau \vert \partial M_i
$$
where $K_i\subset  \partial M_i$ is the preimage of $K$ under the projection map  $\partial M_i\to \abs{\supp \Theta \cap U\cap \partial X}$ defined by $x\mapsto |x|$.
\end{itemize}
\end{thm}
A version of Stokes' theorem holds.

\begin{thm}[Stokes Theorem]
Let $X$ be   an ep-groupoid and let  $\Theta:X\to \Q^+$ be an oriented $n$-dimensional
branched ep-subgroupoid of $X$  whose orbit space  $S=\abs{\supp \Theta}$ is compact. Then,
for every  sc-differential $(n-1)$-form $\omega$ on $X_\infty$,
$$
\mu^{(S,\vartheta)}_{d\omega}(S)=\mu^{(\partial S,\vartheta)}_{\omega}(\partial
S),
$$
or alternatively,
$$\int_{(S, \vartheta)}d\omega =\int_{(\partial S, \vartheta)}\omega.$$
\end{thm}

The construction is compatible with
equivalences between ep-groupoids giving us the following polyfold version.
\begin{thm}
Let $Z$ be a polyfold and $S\subset Z$ be an oriented compact branched suborbifold
defined by the equivalence class $[(X, \beta ,\Theta)]$ and equipped with the weight function $w:S\to \Q^+\cap (0,\infty )$.
For an sc-differential $n$-form $\tau$ on $Z_\infty$ and $K\in {\mathcal L}(S)$,  we  define
$$
\int_{(K,w)}\tau:=\int_{\beta^{-1}(K)}
d\mu_\omega^{(\beta^{-1}(S),\vartheta)}=\mu_\omega^{(\beta^{-1}(S),\vartheta)}(\beta^{-1}(K)),
$$
where  the equivalence class $\tau$ is represented by the triple $(X, \beta ,\omega)$ and where the weight function $\vartheta$  on
$\beta^{-1}(S)=\abs{\supp \Theta}$ is defined by $\vartheta (|x|)=\Theta (x).$
Then    the integral  $\int_{(K,w)}\tau$ is independent of the representative $(X,\beta, \omega)$ in the equivalence class.
Moreover,   if $\tau$ is an sc-differential $(n-1)$-form on $Z_\infty$, then
$$
\int_{(\partial S,w)}\tau =\int_{(S,w)} d\tau.
$$
\end{thm}

In the case of our Gromov-Witten example
consider the subspace $Z_{g,m}$ of $Z$, where $g\geq 0$ and $m\geq 0$ are integers, consisting of all stable curves of arithmetic genus $g$ with $m$ marked points. It has an induced polyfold structure.
For this structure one can show that for every $1\leq i\leq m$  the evaluation map
\begin{eqnarray}\label{xst1}
ev_i:Z_{g,m}\rightarrow Q
\end{eqnarray}
at the $i$-th marked point is  sc-smooth and if $2g+m\geq 3$ then the forgetful map
\begin{eqnarray}\label{xst2}
\sigma:Z_{g,m}\rightarrow \overline{\mathcal M}
\end{eqnarray}
into the Deligne-Mumford space is  sc-smooth as well.
The sc-smoothness of these maps allows us, in particular,
to pull-back differential forms on $Q$ and $\overline{\mathcal M}$ to obtain sc-differential forms on $Z_\infty$. Gromov-Witten invariants can then be obtained by
wedging such pull-backs together and integrating them over
suitable branched suborbifolds obtained as solution sets of the perturbed $\bar{\partial}_J$, which turns out to be a polyfold Fredholm section as discussed below.

\subsection{The Fredholm Package for Polyfolds}

Let $p:W\rightarrow Z$ be a strong polyfold bundle and $f$ an sc-smooth section, which in a suitable model
$P:E\rightarrow X$ is represented by the sc-smooth section functor $F$.
We say that $f$ is {\bf sc-Fredholm} if that is true for $F$. The section $f$ is said to be {\bf proper} provided the solution set of $f=0$ is compact in $Z$ or alternatively the orbit space associated to $F=0$ is compact.

The section $F$ is functorial and a perturbation theory should be functorial as well. Unfortunately functoriality and transversality are two competing concepts. The way out of this dilemma, is a perturbation theory which violates functoriality on the one hand, but tries to
be as compatible as possible with the symmetries. The end result is a  multi-valued perturbation theory incorporating some compatibility with the local symmetries.

We need as before the concept of an auxiliary norm.
\begin{defn}
An {\bf auxiliary norm}\footnote{It is actually not a norm, but comes from an auxiliary norm on a model.} for the strong polyfold bundle $p:W\rightarrow Z$ is a continuous map $N:W_{0,1}\rightarrow [0,\infty)$, so that there exists a model
$P:E\rightarrow X$ and an auxiliary norm $N^\ast$ for $P$ compatible with the morphisms, i.e. if $\phi:h\rightarrow h'$ is a morphism, then $N^\ast(h)=N^\ast(h')$, so that $N\circ \Gamma(|h|) =N^\ast(h)$.
Here $\Gamma:|X|\rightarrow W$ is the identifying homeomorphism.
\end{defn}

If $f$ is sc-Fredholm, the functor $F$ is sc-Fredholm in the same way as discussed before. In particular,  it has the local compactness property.
If $U$ is an open subset of the ep-groupoid $X$, then $|U|$ is an open subset of $|X|$. Since we can find  finitely many (arbitrarily small) open neighborhoods whose orbit spaces cover the compact solution set of a proper section $f$, we can use the local compactness property to construct for given auxiliary norm $N^\ast$ a saturated open neighborhood of $F=0$, so that
$(U,N^\ast)$ controls compactness for functorial $\ssc^+$-sections on the orbit-level, i.e. if $s$ is an $\ssc^+$-section (functorial) with $N^\ast(s(x))<1$ for all $x$ and support in $U$, then the solution set
${\mathcal S}^\ast=\{x\in X\ |\ F(x)=s(x)\}$ has the property that
${\mathcal S}=\gamma(|{\mathcal S}^\ast|)$ is compact in $Z$.

The main point is now to deal with the issue of transversality in this functorial context. In general, i.e. without special structures being present in the problem, we need to use multi-valued perturbations.

 We begin by defining an $\ssc^+$-multisection for a strong bundle over an ep-groupoid. Let us view ${\mathbb Q}^+$ as category only having the identity morphisms.
\begin{defn}
Let $P:E\rightarrow X$ be a strong bundle over the ep-groupoid $X$. An {\bf $\ssc^+$-multisection} is a functor $\Lambda:E\rightarrow {\mathbb Q}^+$, so that for every $x\in X$ there exists an open neighborhood $U$ and finitely many $\ssc^+$-sections $s_1,\ldots ,s_k$ of $E|U$ with associated rational weights $\sigma_i>0$ so that:
$$
\sum_{j=1}^k \sigma_i=1\ \hbox{and}\ \ \Lambda(e)=\sum_{\{i\ |\ s_i(P(e))=e\}} \sigma_i \ \ \hbox{for all} \  \ e\in E|U.
$$
\end{defn}
By convention the sum over the empty set is $0$. An  $\ssc^+$-multisection
for a strong polyfold bundle $p:W\rightarrow Z$ is a map $W\rightarrow {\mathbb Q}^+$ together with an $\ssc^+$-multisection $\Lambda$ for a suitable model $P:E\rightarrow X$, so that
$$
\lambda(\Gamma(|h|))=\Lambda(h).
$$
 More precisely it is, of course, an equivalence class of such pairs
 $(\lambda,\Lambda)$.

If we have two $\ssc^+$-multisections $\lambda_1$ and $\lambda_2$, represented by $\Lambda_1$ and $\Lambda_2$, in general the representatives
$\Lambda_i$ might be defined on different strong bundles, but using a refinement we can pull them back to a common refinement and may assume that they live on the same model. Then we can define the sum
$$
\lambda_1\oplus \lambda_2,
$$
which is represented by
$\Lambda_1\oplus\Lambda_2$, where the latter is defined by the convolution
$$
(\Lambda_1\oplus\Lambda_2)(h)=\sum_{h_1+h_2=h} \Lambda_1(h_1)\Lambda_2(h_2).
$$
We can measure the norm of $\lambda$ with respect to an auxiliary norm
$N$ at a point $z\in Z$ as follows. Take a local representative which comes equipped with $N^\ast$ and $\Lambda$. Pick $x\in X$ so that $\gamma(|x|)=z$
and take a local section structure $(s_i)$ representing $\Lambda$. Then define $N(\lambda,z)$ to be the maximum of the values $N^\ast(s_i(x))$.
One easily verifies that
$$
N(\lambda_1\oplus\lambda_2,z)\leq N(\lambda_1,z)+N(\lambda_2,z).
$$
\begin{prop}
If $f$ is a proper sc-Fredholm section of the strong polyfold bundle $p:W\rightarrow Z$, and $(N,U)$ controls compactness,
then for an $\ssc^+$-multisection $\lambda$ with $N(\lambda,z)<1$ for all $z$ and support in $U$, the solution set
$$
{\mathcal S}(f,\lambda)=\{z\in Z\ |\ \lambda(f(z))>0\}
$$
is compact.
\end{prop}
In applications the main point is to bring the pair $(f,\lambda)$ into a sufficiently general position by slightly perturbing $\lambda$. Note that locally the set ${\mathcal S}(f,\lambda)$ is represented by finitely many problems $F(y)=s_i(y)$, $i=1,\ldots ,k$, where $y\in U(x)$ and each problem carries a rational weight. Solution counts take these weights into consideration. Counting of solutions should, of course, be independent of the particular local section structure.

There is an important fact about the linearized structure at a solution of $\Lambda(F(x))>0$. Namely the finite set of linear sc-Fredholm operators
obtained by linearizing all $F-s_i$ at $x$, for which $F(x)=s_i(x)$,  is independent of the choice of the local section structure, see \cite{HWZ3}. So we can associate to every
point $x$ which solves $\Lambda(F(x))>0$ a collection of sc-Fredholm operators. We write $L(x)$ or $L_{(F,\Lambda)}(x)$ for this collection.

 Now we can generalize certain notions which were previously introduced.
 It makes sense, in order to keep the formalism simple, to formulate everything in terms of $(f,\lambda)$. But, of course, one has to keep in mind the 'overhead', i.e. that everything is defined on the level of models.

\begin{defn}
Let $p:W\rightarrow Z$ be a strong polyfold bundle and $(f,\lambda)$ a pair, where $f$ is an Fredholm section and $\lambda$ an $\ssc^+$-multisection.
Assume that $P:E\rightarrow X$ is a model and $(f,\lambda)$ is represented by $(F,\Lambda)$.
\begin{itemize}
\item[(1)] The pair $(f,\lambda)$ is {\bf transversal pair at the solution $z$} of $\lambda(f(z))>0$ provided that  all sc-Fredholm operators in $L(x)$  are surjective  for a representative $x\in X$. If this holds for all solutions we call $(f,\lambda)$ a {\bf transversal pair}.
 \item[(2)] The pair $(f,\lambda)$ is in {\bf good position at the solution $z$} provided the operators in $L(x)$  for  a representative $x\in X$  are surjective and all their kernel are in good positio to $C_x\subset T_xX$. If this holds for all solutions $z$ we say that $(f,\lambda)$ is in good position.
        \end{itemize}
        \end{defn}

We can view part  (1) of the definition as a special case of (2). If $\partial Z=\emptyset$, then the part (2) specializes to (1).

The first is a global statement about the solution set of a proper Fredholm section satisfying (2)  (or (1)).
\begin{thm}
 Assume that $p:W\rightarrow Z$ is a strong polyfold bundle.  Let $f$ be a proper
Fredholm section and $\lambda$ an $\ssc^+$-multisection so that $(f,\lambda)$ is in good position and the solution set ${\mathcal S}=\{z\in Z\ |\ \lambda(f(z))>0\}$ is compact. Then the pair $({\mathcal S},\lambda\circ f)$ carries in a natural way the structure of a compact branched suborbifold of $Z$ with boundary with corners. Moreover, if $f$ is oriented\footnote{Orientability is again defined via the orientability of an determinant bundle, obtained by linearizing the section. In addition the morphisms are assumed to be compatible with the orientations. Note that the linearized morphisms define maps between the determinants above the points related by a morphism.}, then the branched suborbifold $({\mathcal S},\lambda\circ f)$ is oriented.
\end{thm}

The next result says that we can bring a problem into general position if we allow perturbations by $\ssc^+$-multisec\-tions.

\begin{thm}\label{thmx1}
Let $p:W\rightarrow Z$ be a strong polyfold bundle, where we assume (for simplicity) that $\partial Z=\emptyset$. We assume that $Z$ admits sc-smooth partitions of unity. Let $f$ be a proper
Fredholm section. Assume that $(U,N)$ is a pair controlling compactness for $f$. Then for given $\varepsilon\in (0,1)$ there exists an $\ssc^+$-multisection $\tau$ satisfying $N(\tau)<\varepsilon$ and having support in $U$, so that $(f,\tau)$ is a transversal pair. In particular,  the associated solution set is a branched suborbifold of $Z$.
\end{thm}

Again it is true that given transversal pairs $(f,\lambda_0)$ and $(f,\lambda_1)$, where both perturbations are controlled by $(U_i,N_i)$
there exist a family $\lambda_t$, so that viewing $f$ as a Fredholm section
of $W\rightarrow [0,1]\times Z$ the pair $(f,\lambda_t)$ is transversal, which moreover has the required compactness properties.
Hence we obtain a cobordism between the previous two solution sets.

Let us come back finally to the Gromov-Witten example.
In the Gromov-Witten context $Z$ has sc-smooth partitions of unity which are needed to guarantee a sufficient supply of $\ssc^+$-sections (or corresponding multisections).
Recall the strong bundle $p:W\rightarrow Z$ over the polyfold of stable curves. As before we denote by $Z_{g,m}$ the open subset of $Z$ with its induced polyfold structure consisting of those curves which are connected and have arithmetic genus $g$ and $m$ marked points.

As shown in \cite{HWZ9} (see also \cite{H2}), we have

\begin{thm}
The sc-smooth section $\bar{\partial}_J$ of $p:W\rightarrow Z$ is a polyfold Fredholm section, which is proper on each connected component of $Z$ and is naturally oriented\footnote{The orientation comes from the fact that the linearisations can be convexly homotoped to complex linear Fredholm operators.}.
 \end{thm}

The component-wise properness is, of course, a consequence of Gromov's compactness theorem, \cite{G}.
We can fix a pair $(U,N^\ast)$ controlling compactness, where
$U$ is a suitable open neighborhood of the zero set ${\mathcal S}$ of $\bar{\partial}_J$, and $N^\ast$ an auxiliary norm.
In view of Theorem \ref{thmx1} there are many small perturbations
$\lambda$ of $\bar{\partial}_J$ so that $(\bar{\partial}_J,\lambda)$ is transversal and its solution set
is a naturally oriented branched sub-orbifold, say ${\mathcal S}(\bar{\partial}_J,\lambda)$. In view of the integration theory
we can integrate sc-differential forms over them. In fact not only those coming via pull-backs of the evaluation and forgetful map, but any such form on $Z_\infty$  provided the dimension is right.

\begin{thm}
For every connected component $C$ of $Z$ there exists a uniquely determined linear map $\Phi_C:H^\ast_{dR}(C)\rightarrow {\mathbb R}$ characterized by the property that for a given pair $(U,N^\ast)$ controlling compactness and transversal perturbation $\lambda$ supported in $U$ and $N^\ast(\lambda)<1$ we have
for the solution set $S$ with its natural weight function
$w$
$$
\Phi_C(\tau) =\int_{(S,w)} \tau.
$$
\end{thm}

Using this map and applying it to the wedges of suitably pulled back forms and relating the connected components of $Z$ with
second homology classes of the symplectic manifold $Q$ we can organize the data in various ways, for example as a generating function.

\begin{rem}  There is, of course, a large literature on Gromov-Witten invariants (in the smooth case) and further developments, see \cite{FO,FOOO,LiT,LuT,R1,R2,Tian}. All the methods differ. In some sense, all the approaches had to come up with a fix for the fact that classical Fredholm theory doesn't work.
The theory we describe here seems to come closest to what one would have liked to do in the first place.
\end{rem}


\begin{thebibliography}{99}
\bibitem{AR} A. Adem, J. Leida\ and\ Y. Ruan, Orbifolds and Stringy
Topology, to appear Cambridge University Press.
\bibitem{BSV} V. Borisovich, V. Zvyagin\ and\ V. Sapronov, Nonlinear Fredholm maps and Leray-Schauder degree, Russian
Math. Survey's 32:4 (1977), p 1-54.
\bibitem{BEHWZ} F.~Bourgeois, Y.~Eliashberg, H.~Hofer,
K.~Wysocki and E.~Zehnder,  Compactness Results in Symplectic Field
Theory, {\em Geometry and Topology}, Vol. 7, 2003, pp.799-888.
\bibitem{CRS} K. Cieliebak, I. Mundet i Riera\ and\ D. A. Salamon,
 Equivariant moduli problems, branched manifolds,
 and the Euler class,
Topology {\bf 42} (2003), no.~3, 641--700.
\bibitem{DK} S. Donaldson\ and\ P. Kronheimer, The geometry of
four-manifolds, Oxford Mathematical Monographs. Oxford Science
Publications. The Clarendon Press, Oxford University Press, New
York, 1990.
\bibitem{EGH} Y. Eliashberg, A. Givental\ and\ H. Hofer, Introduction to Symplectic Field Theory,
 Geom. Funct. Anal. {\bf 2000}, Special Volume, Part II, 560--673.
\bibitem{El} H. Eliasson, Geometry of manifolds of maps, J. Differential Geometry {\bf 1}(1967), 169--194.
\bibitem{FH} A. Floer\ and\ H. Hofer, Coherent orientations for periodic orbit problems in symplectic geometry, Math. Z. {\bf 212} (1993), no.~1, 13--38.
    \bibitem{FO} K. Fukaya\ and\ K. Ono,
    Arnold conjecture and Gromov-Witten invariants.
    Topology,Vol. 38 No 5, 1999.pp. 933-1048.
\bibitem{FOOO} K. Fukaya, Y.-G. Oh, H. Ohta and K. Ono,
Lagrangian intersection Floer theory-anomaly and obstruction,
preprint.
\bibitem{GZ} P. Gabriel and M Zisman, Calculus of Fractions and
Homotopy Theory, Ergebnisse Vol. 35, Springer (1967).
\bibitem{G} M.~Gromov, Pseudoholomorphic Curves in
Symplectic Geometry, {\it Inv. Math.} Vol. 82 (1985), 307-347.
\bibitem{hae1} A. Haefliger,  Homotopy and integrability, in Manifolds (Amsterdam, 1970), 133--163, Springer Lecture Notes in Math., 197, 1971.
        \bibitem{hae2} A. Haefliger, Holonomie et classifiants, Astérisque 116 (1984), 70--97.
 \bibitem{hae3} A. Haefliger, Groupoids and foliations, Contemp. Math. 282 (2001), 83--100.
 \bibitem{H1} H. Hofer, A General Fredholm Theory and
Applications, Current Developments in Mathematics, edited by D.
Jerison, B. Mazur, T. Mrowka, W. Schmid, R. Stanley, and S. T. Yau,
International Press, 2006.
\bibitem{H2} H. Hofer, Lectures on Polyfolds and Applications {I}:
Basic Concepts and Illustrations, to appear Courant Lecture Note.
\bibitem{H3} H. Hofer, Lectures on Polyfolds and Applications {II}:
Operations, in preparation.
\bibitem{HWZ-polyfolds1} H. Hofer, K. Wysocki\ and\ E. Zehnder,
Fredholm Theory in Polyfolds {I}: Functional Analytic Methods, book
in preparation.
\bibitem{HWZ-polyfolds2} H. Hofer, K. Wysocki\ and\ E. Zehnder,
Fredholm Theory in Polyfolds {II}: The Polyfolds of Symlectic Field
Theory, book in preparation.
\bibitem{HWZ-DM} H. Hofer, K. Wysocki\ and\ E. Zehnder,
Deligne-Mumford-Type spaces with a View Towards Symplectic Field
Theory, lecture note in preparation.
\bibitem{HWZ1} H. Hofer, K. Wysocki\ and\ E. Zehnder,
A General Fredholm Theory {I}: A Splicing-Based Differential
Geometry,  J. Eur. Math. Soc. (JEMS)  9  (2007),  no. 4, 841--876.
\bibitem{HWZ2} H. Hofer, K. Wysocki\ and\ E. Zehnder,
A General Fredholm Theory {II}: Implicit Function Theorems,
to appear GAFA 2008.
\bibitem{HWZ3} H. Hofer, K. Wysocki\ and\ E. Zehnder,
A General Fredholm Theory {III}: Fredholm Functors and Polyfolds, preprint.
\bibitem{HWZ4} H. Hofer, K. Wysocki\ and\ E. Zehnder,
A General Fredholm Theory {IV}: Operations , paper in preparation.
\bibitem{HWZ5} H. Hofer, K. Wysocki\ and\ E.
Zehnder, Integration Theory for  zero sets of polyfold Fredholm
sections,
preprint.
\bibitem{HWZ6} H. Hofer, K. Wysocki\ and\ E. Zehnder,
Applications of Polyfold Theory {I}: Gromov-Witten Theory, paper in
preparation.
\bibitem{HWZ7} H. Hofer, K. Wysocki\ and\ E. Zehnder,
Applications of Polyfold Theory {II}: The Polyfolds of Symplectic
Field Theory, paper in preparation.
\bibitem{HWZ8} H. Hofer, K. Wysocki\ and\ E. Zehnder,
Connections and Determinant Bundles for Polyfold Fredholm Operators,
paper in preparation.
\bibitem{HWZ9} H. Hofer, K. Wysocki\ and\ E. Zehnder,
The Fredholm Property of the Nonlinear Cauchy-Riemann operator in
the Polyfold Set-up for Symplectic Field Theory, paper in
preparation.
\bibitem{Joyce} D. Joyce, Kuranishi bordism and Kuranishi homology,
preprint arxiv 0707.3572v1.
\bibitem{LA} S. Lang, {\it
Introduction to differentiable manifolds}, Second edition, Springer,
New York, 2002.
\bibitem{LiT} J. Li and G. Tian, Virtual moduli cycles and Gromov-Witten
invariants of general symplectic manifolds, in: Topics in symplectic
4-manifolds (Irvine, CA, 1996), 47–83, Int. Press (1998).
\bibitem{LuT} G. Lu and G. Tian, Constructing virtual Euler cycles and classes,  Int. Math. Res. Surv. IMRS  2007,  2008 in electronic version, Art. ID rym001, 220 pp.
\bibitem{Man} Y. Manin, {\it Frobenius Manifolds, Quantum Cohomology, and
Moduli Spaces}, AMS Colloquium Publications, Volume 47.
\bibitem{Mc} D. McDuff, Groupoids, Branched Manifolds and
Multisection, J. Symplectic Geom. 4, 259-315 (2006).
\bibitem{MS} D. McDuff and D. Salamon, {\it Introduction
 to symplectic topology}, 2nd edition, Oxford University Press, 1998.
\bibitem{MS2} D. McDuff and D. Salamon, {\it J-holomorphic curves and
symplectic topology} , Colloquium Publications, vol. 52, Amer. Math.
Soc., Providence, RI, 2004, xii+669 pp.
\bibitem{Mj} I. Moerdijk, Orbifolds as Groupoids: An Introduction, Contemp. Math. 310, 205-222 (2002).
\bibitem{MM} I. Moerdijk and J. Mr\v cun, {\it Introduction to
Foliation and Lie Groupoids}, Cambridge studies in advanced
mathematics, Vol. 91, 2003.
\bibitem{R1} Y. Ruan, Topological sigma model and Donaldson type
invariants in Gromov theory, Duke Math. J. 83, 1996, 451-500.
\bibitem{R2} Y. Ruan, Symplectic Topology on Algebraic 3-folds, J.
Diff. Geom. 39, 1994, 215-227.
\bibitem{smale} S. Smale,
An infinite dimensional version of Sard's theorem.
Amer. J. Math. 87 (1965), 861--866.
\bibitem{Tian} G. Tian, Quantum cohomology and its associativity,
Current Developments in Mathematics, 1995, p. 361-397, International
Press
\bibitem{Tr} H. Triebel, {\it Interpolation theory, function spaces, differential operators}, North-Holland, Amsterdam, 1978.


\end{thebibliography}
\end{document}